\documentclass[reqno]{amsart}

\usepackage{amsmath, amssymb}
\usepackage{amsthm, stmaryrd}
\usepackage{mathrsfs}

\usepackage[top=40truemm,bottom=40truemm,left=30truemm,right=30truemm]{geometry}

\usepackage{hyperref}

\usepackage{url}

\makeatletter
	\@addtoreset{equation}{section}

\theoremstyle{plain}
\newtheorem{theorem}{Theorem}[section]
\newtheorem{proposition}[theorem]{Proposition}
\newtheorem{lemma}[theorem]{Lemma}
\newtheorem{corollary}[theorem]{Corollary}

\theoremstyle{definition}
\newtheorem{definition}[theorem]{Definition}
\newtheorem{remark}[theorem]{Remark}
\newtheorem{example}{Example}[section]

\allowdisplaybreaks[4]

\date{}

\author{Takumu Ooi}
\address{Department of Mathematics,
Faculty of Science and Technology, Tokyo University of Science,
Noda, Chiba, 278-8510, Japan}
\email{ooitaku@rs.tus.ac.jp}

\title[Analytic criteria for criticality of subordinated Schr\"{o}dinger op.]{Analytic criteria for criticality of subordinated Schr\"{o}dinger operators}

\begin{document}
\maketitle

\begin{abstract}
For a subordinated Schr\"{o}dinger operator, we introduce a quantity that plays a role analogous to the bottom of the spectrum of a Schr\"{o}dinger operator, and establish analytic criteria for subcriticality, criticality and supercriticality of a subordinated Schr\"{o}dinger operator.
\end{abstract}

\section{Introduction}
For Schr\"{o}dinger operators, which are obtained by perturbing non-positive self-adjoint operators with signed measures, their subcriticality, criticality and supercriticality have been extensively studied within the framework of Schr\"{o}dinger form theory based on Dirichlet form theory (see \cite{T14, TU23} for example). More recently, \cite{OS26} introduced probabilistic notions of subcriticality, criticality and supercriticality for subordinated Schr\"{o}dinger operators such as \((-\Delta + V)^{\alpha/2}\) (see also \cite{S26}). They also established a connection between the subcriticality of these operators and the boundedness of solutions to certain wave-type equations. However, determining subcriticality through its probabilistic characterization is generally a challenging task. In this paper, we introduce an analytic quantity for subordinated Schr\"{o}dinger operators that plays a role analogous to the bottom of the spectrum in the theory of Schr\"{o}dinger operators. Using this quantity, we establish analytic characterizations of subcriticality and criticality for subordinated Schr\"{o}dinger operators.

Let \((\mathcal{E}, \mathcal{D}(\mathcal{E}))\) be an irreducible regular Dirichlet form and \(\mu\) be a signed smooth Radon measure such that \((\mathcal{E}^{\mu}, \mathcal{D}(\mathcal{E})\cap C_c(E))\) is a non-negative definite closable symmetric form. Then its closure \((\mathcal{E}^\mu, D(\mathcal{E}^\mu))\) corresponds to a non-positive self-adjoint operator \(\mathcal{L}^\mu\). For a Bernstein function \(\Phi\), we consider a subordinated Schr\"{o}dinger operator \(-\Phi(-\mathcal{L}^\mu)\). Let \((\mathcal{E}^{\mu,\Phi}, D(\mathcal{E}^{\mu,\Phi}))\) be a closed bilinear form associated with \(-\Phi(-\mathcal{L}^\mu)\). Under the condition making \((\mathcal{E}^{\mu,\Phi}, D(\mathcal{E}^{\mu,\Phi}))\)  irreducible, it is defined and proved in \cite{OS26} that \(-\Phi(-\mathcal{L}^\mu)\) is subcritical (resp. critical) if a Doob's \(h\)-transform of \((\mathcal{E}^{\mu,\Phi}, D(\mathcal{E}^{\mu,\Phi}))\) is transient (resp. recurrent) for some (equivalently any) strictly positive excessive function \(h\). 
For the case of Schr\"{o}dinger operators, it is known that, the spectral bottom \(\lambda(\mu)\) of the time-changed process of \(\mathcal{E}^{\mu_+}\) by \(\mu_-\) characterize the subcriticality and criticality. Here \(\mu_+\) (resp. \(\mu_-\)) is a positive (resp. negative) part of \(\mu\). More precisely, if \(\lambda(\mu)>1\) then a Schr\"{o}dinger operator \(\mathcal{L}^\mu\) is subcritical (\cite{TU23}). When \(\mu_+\) is a local Kato class measure and \(\mu_-\) is a Green-tight measure for \(\mathcal{E}^\mu_+ \), if \(\lambda(\mu)=1\) then a Schr\"{o}dinger operator \(\mathcal{L}^\mu\) is critical (\cite{T14}). In this paper, for a complete Bernstein function \(\Phi\), we introduce \(\lambda(\mu, \Phi)\) by
\[\lambda(\mu, \Phi):=\inf_{f\in D(\mathcal{E}) \cap C_c(E)} \frac{\mathcal{E}^{\mu_+, \Phi}(f,f)}{\mathcal{E}^{\mu_+, \Phi}(f,f) - \mathcal{E}^{\mu, \Phi}(f,f)}\]
and we prove that the following analytic criteria:
\begin{theorem}[Theorem \ref{thm:subcri}]
If \(\lambda(\mu, \Phi)>1\) then a subordinated Schr\"{o}dinger operator \(-\Phi(-\mathcal{L}^\mu)\) is subcritical.
\end{theorem}

\begin{theorem}[Theorem \ref{thm:cri}]
Assume that \((\mathcal{E}, D(\mathcal{E}))\) satisfies the strong Feller property, \(\mu_+\) is a local Kato class measure, \(\mu_-\) is a Green-tight measure for \(\mathcal{E}^{\mu_+} \) and the condition {\rm (UC)}. If \(\lambda(\mu, \Phi)=1\) then a subordinated Schr\"{o}dinger operator \(-\Phi(-\mathcal{L}^\mu)\) is critical.
\end{theorem}
The condition {\rm (UC)} will be defined in Section \ref{sec:main}, and we also provide the sufficient conditions for {\rm (UC)}. By using these criteria, we can distinguish criticality and subcriticality for subordinated Schr\"{o}dinger operators.

The organization of this paper is as follows. In Section \ref{sec:pre}, we review previous work \cite{T14, TU23} for Schr\"{o}dinger forms and \cite{OS26} for subordinated Schr\"{o}dinger forms. In Section \ref{sec:main}, we introduce \(\lambda(\mu, \Phi)\) for a subordinated Schr\"{o}dinger form and consider properties for \(\lambda(\mu, \Phi)\). We also prove the main results Theorem 1.1, 1.2, criteria for criticality and subcriticality for subordinated Schr\"{o}dinger forms. Section \ref{sec:Example} gives examples. In particular, for a subordinated Schr\"{o}dinger form corresponding to the critical Hady inequality for the Laplacian and the fractional Laplacian on \(\mathbb{R}^d\), we classifies \(\lambda(\mu, \Phi)\) for all complete Bernstein functions. Also we consider the case of bounded integrable potentials. Appendix \ref{append} contains preliminaries on Dirichlet form theory and Markov
processes.

\section{Preliminaries}\label{sec:pre}
Let \((\mathcal{E}, \mathcal{D}(\mathcal{E}))\) be an irreducible regular Dirichlet form and \(\mu\) a signed smooth Radon measure such that \((\mathcal{E}^{\mu}, \mathcal{D}(\mathcal{E})\cap C_c(E))\) is non-negative definite closable symmetric form. For a signed smooth Radon measure \(\mu\), we denote by \(\mu_+\) (resp. \(\mu_-\)) a positive (resp. negative) part of \(\mu\), and, we write \(\mu \in \mathcal{S}_+-\mathcal{S}_-\) when \(\mu_+ \in \mathcal{S}_+\) and \(\mu_- \in \mathcal{S}_-\) for classes \(\mathcal{S}_+, \mathcal{S}_- \subset \mathcal{S}\).

Denote by \((\mathcal{E}^{\mu}, \mathcal{D}(\mathcal{E}^\mu))\) its Schr\"{o}dinger form 
, \(\mathcal{L}^{\mu}=\mathcal{L}-\mu\) its associated non-positive self-adjoint operator, and \(\{T_t^\mu\}_t\) its strongly continuous contraction semigroup on \(L^2(E;m)\). We define the resolvent \(\{R_\alpha^\mu\}_{\alpha>0}\) on \(L^2(E;m)\) by 
\[R_\alpha^\mu f := \int_0^\infty e^{-\alpha t}T_t^\mu f\,dt,\] 
and it holds that \(\mathcal{E}_\alpha^\mu(R_\alpha^\mu f, g)=\langle f,g \rangle_m\) where \(\mathcal{E}_\alpha^\mu(f,g):=\mathcal{E}^\mu(f,g) + \alpha \langle f, g\rangle_m.\)

Without loss of generality, we assume that \((\mathcal{E}^{\mu_+}, D(\mathcal{E}^{\mu_+}))\) is transient and \(\mu_-\) is non-zero. Indeed, if \((\mathcal{E}^{\mu_+}, D(\mathcal{E}^{\mu_+}))\) is not transient, then, by irreducibility, this is recurrent and so \(1 \in D_e(\mathcal{E}^{\mu_+})\), where \(D_e(\mathcal{E}^{\mu_+})\) is an extended Dirichlet space. Since \((\mathcal{E}^{\mu}, \mathcal{D}(\mathcal{E})\cap C_c(E))\) is non-negative definite, we have \(0=\mathcal{E}^{\mu_+}(1,1)\geq \mu_-(E)\). Hence \((\mathcal{E}^\mu, \mathcal{D}(\mathcal{E}^\mu))= (\mathcal{E}^{\mu_+}, \mathcal{D}(\mathcal{E}^{\mu_+}))\) is a regular Dirichlet form, and we can apply the classical theory of Dirichlet forms. 

Let \(S\) be a conservative subordinator, which is a non-decreasing L\'{e}vy process on \([0,\infty)\). For the probability distribution  \(\eta_t\) of \(S_t\), denote by \(\Phi\) a Laplace exponent. More precisely, \(\int_0^\infty e^{-\lambda s}\,\eta_t(ds) = e^{-\lambda \Phi(t)}\) holds. By the L\'{e}vy-Khintchine formula (see, for example, \cite[Theorem 8.1]{S99}), we have
\begin{eqnarray*}
\Phi(\lambda)= b\lambda + \int_0^\infty (1-e^{-\lambda s})\, \nu(ds),
\end{eqnarray*}
where \(b \geq 0\) is a constant called the drift coefficient and \(\nu\) is a Radon measure on \((0, \infty)\) satisfying \(\int_0^\infty (1\wedge s)\,d\nu(s) < \infty\), called the  L\'{e}vy  measure. \(\Phi\) is called a {\it Bernstein function.} A Bernstein function is a smooth function satisfying \((-1)^{n-1} \frac{d^n}{d\lambda^n}\Phi(\lambda) \geq 0\) for any \(n\in \mathbb{N}\) and \(\lambda\geq 0.\) See \cite{SSV12} for details on Bernstein functions.

We can define a strongly continuous contraction semigroup \(\{T_t^{\mu, \Phi}\}_t\) on \(L^2(E;m)\) by
\[T_t^{\mu, \Phi}f :=\int_0^\infty T_s^{\mu}f\, \eta_t(ds) \]
for \(f\in L^2(E;m)\). The generator \(\mathcal{L}^{\mu, \Phi}\) of this semigroup is represented by 
\[\mathcal{L}^{\mu, \Phi} = -\Phi(-\mathcal{L}^{\mu}) = b\mathcal{L}^{\mu} - \int_0^\infty (I-T_s^{\mu})\, \nu(ds)\]
and we have \(\mathcal{D}(\mathcal{L}^{\mu, \Phi}) \supset \mathcal{D}(\mathcal{L}^{\mu}).\) We call \(\mathcal{L}^{\mu, \Phi}\) a subordinated Schr\"{o}dinger operator, and its closed form \((\mathcal{E}^{\mu, \Phi}, \mathcal{D}(\mathcal{E}^{\mu, \Phi}))\) on \(L^2(E;m)\), which is given by \(\mathcal{E}^{\mu, \Phi}(f,g):=\langle \sqrt{-\mathcal{L}^{\mu, \Phi}}f, \sqrt{-\mathcal{L}^{\mu, \Phi}}g \rangle_m \) and \(\mathcal{D}(\mathcal{E}^{\mu, \Phi}) = \mathcal{D}(\sqrt{-\mathcal{L}^{\mu, \Phi}})\), {\it a subordinated Schr\"{o}dinger form}. We set \(\mathcal{E}_\alpha ^{\mu, \Phi}(f,g):= \mathcal{E}^{\mu, \Phi}(f,g) +\alpha \langle f,g \rangle_m \).

By \cite[Theorem 3.2]{OS26} (See also \cite[Theorem 2.1]{O02}), it holds that \(\mathcal{D}(\mathcal{E}^{\mu}) \subset \mathcal{D}(\mathcal{E}^{\mu, \Phi})\) and the subordinated Schr\"{o}dinger form \((\mathcal{E}^{\mu, \Phi}, \mathcal{D}(\mathcal{E}^{\mu, \Phi}))\) on \(L^2(E;m)\) admits the representation
\[\mathcal{E}^{\mu, \Phi}(f,g)=b\, \mathcal{E}^{\mu}(f,g) + \int_0^\infty \langle f-T_s^\mu f, g \rangle_m d\nu(s)\]
for \(f,g \in \mathcal{D}(\mathcal{E}^{\mu})\).

We define the extended space \(\mathcal{D}_e(\mathcal{E}^{\mu, \Phi})\) by the set of all \(m\)-measurable functions \(f\) satisfying \(|f|<\infty\) \(m\)-almost everywhere and possessing an approximating sequence \(\{f_n\}_n \subset \mathcal{D}(\mathcal{E}^{\mu, \Phi})\) such that \(\mathcal{E}^{\mu, \Phi}(f_n-f_m, f_n-f_m)\to 0\) as \(m,n \to \infty\) and \(f_n \to f\) \(m\)-almost everywhere, and we can define \(\mathcal{E}^{\mu, \Phi}(f,f):=\lim_{n\to \infty} \mathcal{E}^{\mu, \Phi}(f_n,f_n)\). We note that the limit  \(\mathcal{E}^{\mu, \Phi}(f,f)\) is independent of the choice of an approximating sequence \(\{f_n\}_n\). See \cite[Section 3]{OS26} for details.

\begin{definition}[{\cite[Definition 3.3]{OS26}}]\, 
\begin{enumerate}
\item  A subordinated Schr\"{o}dinger form \((\mathcal{E}^{\mu, \Phi}, \mathcal{D}(\mathcal{E}^{\mu, \Phi}))\) is subcritical if there exists a strictly positive bounded function \(g \in L^1(E;m)\) such that, for any \(f\in \mathcal{D}_e(\mathcal{E}^{\mu, \Phi})\),
\begin{equation} \int_E |f|g \,dm \leq \sqrt{\mathcal{E}^{\mu, \Phi}(f,f)}. \label{eq:sub_subcri} \end{equation}
\item A subordinated Schr\"{o}dinger form \((\mathcal{E}^{\mu, \Phi}, \mathcal{D}(\mathcal{E}^{\mu, \Phi}))\) is critical if there exists a strictly positive function \(h\in  \mathcal{D}_e(\mathcal{E}^{\mu, \Phi})\) satisfying \(\mathcal{E}^{\mu, \Phi}(h,h)=0.\)
\item  A subordinated Schr\"{o}dinger form \((\mathcal{E}^{\mu, \Phi}, \mathcal{D}(\mathcal{E}^{\mu, \Phi}))\) is supercritical if neither \((1)\) nor \((2)\) holds.
\end{enumerate}
Moreover, we say that \(\mathcal{L}^{\mu, \Phi}\) is subcritical (resp. critical, supercritical) if \((\mathcal{E}^{\mu, \Phi}, \mathcal{D}(\mathcal{E}^{\mu, \Phi}))\) is subcritical (resp. critical, supercritical).
\end{definition}

We define the set \(\mathcal{H}^{\mu, \Phi}_+\) of all strictly positive excessive functions by
 \[\mathcal{H}^{\mu, \Phi}_+:=\{h : 0<h<\infty,\  T_t^{\mu, \Phi}h \leq h \text{\ for\ any\ }t>0, m\text{-almost everywhere}\}.\]
Then it holds that \(\mathcal{H}_+^\mu \subset \mathcal{H}^{\mu, \Phi}_+ \). For \(h \in \mathcal{H}^{\mu, \Phi}_+\), we define a Doob's \(h\)-transform \((\mathcal{E}^{\mu,\Phi, h}, \mathcal{D}(\mathcal{E}^{\mu,\Phi,h}))\) by
\begin{eqnarray*}
\mathcal{E}^{\mu,\Phi,h}(f,g)&:=& \mathcal{E}^{\mu, \Phi}(fh,gh)\\
\mathcal{D}(\mathcal{E}^{\mu,\Phi,h})&:=& \{f\in L^2(E;h^2m) : fh\in \mathcal{D}(\mathcal{E}^{\mu, \Phi})\}.
\end{eqnarray*}

We consider the condition (IB) as follows.

(IB): A Bernstein function \(\Phi\) satisfies  either \(b>0\) or \(\inf{{\rm supp}(\nu)}=0\).\\

Under the condition (IB), \((\mathcal{E}^{\mu,\Phi}, \mathcal{D}(\mathcal{E}^{\mu,\Phi}))\) is irreducible (\cite[Proposition 3.4]{OS26}) and the following characterizations of subcriticality, criticality, and supercriticality hold. 

\begin{theorem}[{\cite[Theorem 3.5]{OS26}}]\label{subodinated_subcri}
Under the condition \rm{(IB)}, the following are equivalent.
\begin{enumerate}
\item A subordinated Schr\"{o}dinger form \((\mathcal{E}^{\mu, \Phi}, \mathcal{D}(\mathcal{E}^{\mu, \Phi}))\) is  subcritical.
\item \(\mathcal{H}_+^{\mu, \Phi}\) is not empty and \((\mathcal{E}^{\mu,\Phi, h}, \mathcal{D}(\mathcal{E}^{\mu, \Phi,h}))\) on \(L^2(E;h^2m)\) is transient for some \(h\in \mathcal{H}_+^{\mu, \Phi}\).
\item \(\mathcal{H}_+^{\mu, \Phi}\) is not empty and \((\mathcal{E}^{\mu,\Phi, h}, \mathcal{D}(\mathcal{E}^{\mu, \Phi,h}))\) on \(L^2(E;h^2m)\) is transient for any \(h\in \mathcal{H}_+^{\mu, \Phi}\).
\end{enumerate}
\end{theorem}

\begin{theorem}[{\cite[Theorem 3.6]{OS26}}]\label{subodinated_cri}
Under the condition \rm{(IB)}, the following are equivalent.
\begin{enumerate}
\item A subordinated Schr\"{o}dinger form \((\mathcal{E}^{\mu, \Phi}, \mathcal{D}(\mathcal{E}^{\mu, \Phi}))\) is  critical.
\item \(\mathcal{H}_+^{\mu, \Phi}\) is not empty and \((\mathcal{E}^{\mu,\Phi, h}, \mathcal{D}(\mathcal{E}^{\mu, \Phi,h}))\) on \(L^2(E;h^2m)\) is recurrent for some \(h\in \mathcal{H}_+^{\mu, \Phi}\).
\item \(\mathcal{H}_+^{\mu, \Phi}\) is not empty and \((\mathcal{E}^{\mu,\Phi, h}, \mathcal{D}(\mathcal{E}^{\mu, \Phi,h}))\) on \(L^2(E;h^2m)\) is recurrent for any \(h\in \mathcal{H}_+^{\mu, \Phi}\).
\end{enumerate}
\end{theorem}

\begin{theorem}[{\cite[Theorem 3.7]{OS26}}]\label{subodinated_supercri}
Under the condition \rm{(IB)}, the following are equivalent.
\begin{enumerate}
\item A subordinated Schr\"{o}dinger form \((\mathcal{E}^{\mu, \Phi}, \mathcal{D}(\mathcal{E}^{\mu, \Phi}))\) is  supercritical.
\item \(\mathcal{H}_+^{\mu, \Phi}\) is empty.
\end{enumerate}
\end{theorem}

These results provide probabilistic characterizations of notions of criticality. We will consider an analytic characterization of notions of criticality from the viewpoint of spectral theory. Before considering an analytic characterization for subordinated Schr\"{o}dinger operators, we recall a characterization for Schr\"{o}dinger operators due to Takeda and Uemura \cite{T14, TU23}. We define the bottom of the spectrum \(\lambda(\mu)\) for the time-changed process of a process associated with \((\mathcal{E}^{\mu_+}, \mathcal{D}(\mathcal{E}^{\mu_+}))\) by \(\mu_-\) as follows.
\[\lambda(\mu):= \inf_{0\neq f\in \mathcal{D}(\mathcal{E})\cap C_c(E)}\frac{\mathcal{E}^{\mu_+}(f,f)}{\int_E |f|^2\, d\mu_-}.\]

We note that \(\mathcal{D}(\mathcal{E}^{\mu_+}) = \mathcal{D}(\mathcal{E}) \cap L^2(E;\mu_+)\) by \cite[Section 5.1]{CF12}.
\begin{theorem}[{\cite[Theorem 3.5]{TU23}}]\label{Theorem3.5_TU23}
Suppose that, for each compact set \(K\),
\[\sup_{f\in \mathcal{D}_e(\mathcal{E}^{\mu_+})} \frac{\left|\int_K f\,dm \right|}{\sqrt{\mathcal{E}^{\mu_+}(f,f)}} < \infty. \]
If \(\lambda(\mu)>1\), then \((\mathcal{E}^\mu, \mathcal{D}(\mathcal{E}^\mu))\) is subcritical.
\end{theorem}

We consider the following the strong Feller property.
\begin{enumerate}
\item[(SF)] : For each \(t>0\), \(T_t (\mathcal{B}_b(E)) \subset C_b(E)\), where \(\mathcal{B}_b(E)\) is the set of all bounded Borel measurable functions on \(E\) and \(C_b(E)\) is the set of all bounded continuous functions on \(E\).
\end{enumerate}
We say that the strong Feller property {\rm (SF)} holds for \((\mathcal{E}, D(\mathcal{E}))\) if the semigroup associated with  \((\mathcal{E}, D(\mathcal{E}))\) satisfies {\rm (SF)}.

By \((\)SF\()\) and the symmetry of the Dirichlet form, the following condition holds.
\begin{enumerate}
\item[(AC)] Absolute continuity condition: There exists a jointly measurable function \(p(t,x,y)\) on \((0,\infty)\times E \times E\) such that \(p(t,x,y)=p(t,y,x)\) and \(P_tf(x)=\int_E p(t,x,y)f(y)dm(y)\).
\end{enumerate}

We define\[r_{\alpha}(x,y) := \int_0^{\infty} e^{-\alpha t} p(t,x,y)\,dt\]
for \(\alpha>0\). Then, for an \(\alpha\)-order resolvent \(R_{\alpha}\), it holds that \(R_{\alpha}f(x)=\int_E r_{\alpha}(x,y)f(y)\,dm(y)\). We note that \(T_tf=P_tf\) and \(G_{\alpha}f=R_{\alpha}f\) hold \(m\)-almost everywhere for \(f\in L^2(E;m)\).

We remark that {\rm (SF)} for  \((\mathcal{E}, D(\mathcal{E}))\) follows from {\rm (SF)} for  \((\mathcal{E}^{\mu_+}, D(\mathcal{E}^{\mu_+}))\). Indeed, by the Dini theorem and the strong Feller property for \((\mathcal{E}^{\mu_+}, D(\mathcal{E}^{\mu_+}))\), we have
\begin{eqnarray*}
\varlimsup_{s \searrow 0} \sup_{x\in K}|T_tf(x) - T_s^{\mu_+}T_{t-s}f(x)| &\leq & \varlimsup_{s \searrow 0} \sup_{x\in K} \mathbb E_x\left[(1-e^{-A_s^{\mu_+}})T_{t-s}f (X_s);\,s<\zeta \right]\\
&\leq & \|f\|_\infty\, \varlimsup_{s \searrow 0} \sup_{x\in K}  |T_s1(x)-T_s^{\mu_+}1(x)| \\
&\leq & \|f\|_\infty\, \varlimsup_{s \searrow 0} \sup_{x\in K}  |1-T_s^{\mu_+}1(x)|=0
\end{eqnarray*}
for a compact set \(K\) and a bounded Borel measurable function \(f\), where \(X\) is a Hunt process associated with \((\mathcal{E}, D(\mathcal{E}))\), \(\zeta\) is a lifetime of \(X\), and \(A^{\mu_+}\) is a PCAF corresponding to \(\mu_+\). Hence the strong Feller property for \((\mathcal{E}, D(\mathcal{E}))\) holds.

\begin{definition}\label{defKato}
Suppose that the condition (AC) holds and we set \(R_{\alpha}\mu(x) := \int_E r_{\alpha}(x,y)\, d\mu(y) \leq \infty\) for \(\mu \in \mathcal{S}\). Then we define the following classes of smooth measures.
\begin{enumerate}
\item A smooth measure \(\mu \in \mathcal{S}\) is called a {\it Kato class measure} if \(\|R_{\alpha}\mu\|_\infty \to 0\) as \(\alpha \to \infty.\) Denote by \(\mathcal{K}:=\mathcal{K}(\mathcal{E})\) the set of all Kato class measures.

\item A smooth measure \(\mu \in \mathcal{S}\) is called a {\it local Kato class measure} if, for any compact set \(K\), the restriction \(1_{K}\, \mu\) of \(\mu\) to \(K\) is a Kato class measure. Denote by \(\mathcal{K}_{loc}:=\mathcal{K}_{loc}(\mathcal{E})\) the set of all local Kato class measures.

\item Suppose that \((\mathcal{E}, \mathcal{D}(\mathcal{E}))\) is transient. A smooth measure \(\mu \in \mathcal{K}\) satisfying \(\|R\mu\|_\infty < \infty\) is called a {\it Green-tight measure} with respect to \((\mathcal{E}, \mathcal{D}(\mathcal{E}))\), if for any \(\varepsilon>0,\) there exists a compact set \(K\) such that \(\|R(1_{K^c}\mu)\|_\infty <\varepsilon\) holds. Denote by \(\mathcal{K}_{\infty}(\mathcal{E})\) the set of all Green-tight measures with respect to \((\mathcal{E}, \mathcal{D}(\mathcal{E}))\).
\end{enumerate}
\end{definition}

For \(\mu_- \in \mathcal{K}(\mathcal{E}^{\mu_+})\), by the Stollmann-Voigt inequality (See \cite{SV96} and \cite[Section 3]{KOTU26} for example), the closability of \((\mathcal{E}^\mu, \mathcal{D} \cap C_c(E))\) automatically holds and we have the identity \(\mathcal{D}(\mathcal{E}^\mu) =\mathcal{D}(\mathcal{E}^{\mu_+})\).

\begin{theorem}[{\cite[Theorem 5.19]{T14}}] \label{Takedaeigen}
Suppose that the strong Feller property {\rm (SF)} holds for \((\mathcal{E}^{\mu_+}, D(\mathcal{E}^{\mu_+}))\). For a signed smooth Radon measure \(\mu \in \mathcal{K}_{loc}(\mathcal{E})-\mathcal{K}_{\infty}(\mathcal{E}^{\mu_+})\), the following holds.
\begin{enumerate}
\item \((\mathcal{E}^{\mu}, \mathcal{D}(\mathcal{E}^{\mu}))\) is subcritical if and only if \(\lambda(\mu)>1\),
\item \((\mathcal{E}^{\mu}, \mathcal{D}(\mathcal{E}^{\mu}))\) is critical if and only if \(\lambda(\mu)=1\),
\item \((\mathcal{E}^{\mu}, \mathcal{D}(\mathcal{E}^{\mu}))\) is supercritical if and only if \(\lambda(\mu)<1\).
\end{enumerate}
\end{theorem}

\section{Analytic criteria}\label{sec:main}
In this section, we consider analytic criteria for subcriticality and criticality of subordinated Schr\"{o}dinger operators from a spectral viewpoint. As in the previous section, we consider an irreducible regular Dirichlet form \((\mathcal{E}, \mathcal{D}(\mathcal{E}))\) on \(L^2(E;m)\) and a signed smooth Radon measure \(\mu\)  such that \((\mathcal{E}^{\mu}, \mathcal{D}(\mathcal{E})\cap C_c(E))\) is a non-negative definite closable symmetric form. Denote by \((\mathcal{E}^{\mu}, \mathcal{D}(\mathcal{E}^\mu))\) its Schr\"{o}dinger form and \(\mathcal{L}^{\mu}=\mathcal{L}-\mu\) its associated non-positive self-adjoint operator. 

In this section, we consider the following type of Bernstein functions.
\begin{definition}[{\cite[Definition 1.3 and 6.1]{SSV12}}]
A Bernstein function \(\Phi\) is called a complete Bernstein function if there exists a \(C^\infty(0,\infty)\)-function \(\rho\) such that \(\nu(ds)= \rho(s) \,ds\) holds for the L\`{e}vy measure \(\nu\), and \((-1)^n \frac{d^n}{ds^n}\rho(s)\geq 0\) for any \(n \in \mathbb{N} \cup \{0\}\) and \(s >0\). 
\end{definition}

By Bernstein's theorem, the density of the Lévy measure of a complete Bernstein function admits (\cite[Theorem 1.4]{SSV12}), that is, there exists a unique measure \(\widetilde{\nu}\) on \((0,\infty)\) such that, for all \(\alpha>0\), 
\[\rho(s) = \int_0^\infty e^{-\alpha s} d\widetilde{\nu}(\alpha).\]
Since \(\int_0^\infty (s\wedge 1)\,d\nu(s) < \infty\), we have
\begin{equation}
\int_0^\infty \frac{1-e^{-\alpha}}{\alpha^2} d\widetilde{\nu}(\alpha) < \infty \quad \big(\text{ or\ equivalently,} \int_0^\infty \frac{1}{\alpha (1+\alpha)} d\widetilde{\nu}(\alpha) < \infty\big). \label{eq:Levytilde}
\end{equation}

A complete Bernstein function \(\Phi\) with \(\Phi \not =0\) satisfies the condition (IB). Indeed, if \(b=0\), \(\inf {\rm supp}(\nu)=0\) follows from \(\rho(s)>0\) for \(s>0\).

We introduce the quantity \(\lambda(\mu, \Phi)\) for a subordinated Schr\"{o}dinger operator.
\begin{definition}
For a complete Bernstein function \(\Phi\), we define
\[\lambda(\mu, \Phi):=\inf_{0\neq f\in D(\mathcal{E}) \cap C_c(E)} \frac{\mathcal{E}^{\mu_+, \Phi}(f,f)}{\mathcal{E}^{\mu_+, \Phi}(f,f) - \mathcal{E}^{\mu, \Phi}(f,f)}.\]
\end{definition}

\begin{remark}
(1) By \cite[Theorem 3.2]{OS26}, \(D(\mathcal{E})\cap C_c\) is dense in  \(D(\mathcal{E}^{\mu,\Phi})\) with respect to \(\mathcal{E}_1^{\mu,\Phi}\) and \(D(\mathcal{E}^{\mu_+,\Phi})\) with respect to \(\mathcal{E}_1^{\mu_+,\Phi}\).

(2) Since \(\mathcal{E}^{\mu}(f,f) \leq \mathcal{E}^{\mu_+}(f,f)\) for any \(f \in D(\mathcal{E}) \cap C_c(E)\) and \(\Phi\) is a complete Bernstein function, by \cite[Theorem 12.12]{SSV12}, we have \(\mathcal{E}^{\mu, \Phi}(f,f) \leq \mathcal{E}^{\mu_+, \Phi}(f,f)\). Hence it holds that \(1 \leq \lambda(\mu, \Phi).\)

(3) For a Bernstein function that is not complete,  \(\mathcal{E}^{\mu, \Phi}(f,f) \leq \mathcal{E}^{\mu_+, \Phi}(f,f)\) does not hold in general. For example, let \(E:=\{0,1\}\), \(m\) be a counting measure, \(\mathcal{E}(f,g):=2(f(1)-f(0))(g(1)-g(0))+3f(1)g(1)\), \(D(\mathcal{E})=L^2(E;m) \simeq \mathbb{R}^2\) and \(\mu:=-\mu_-:=-3\delta_1\). Then, \((\mathcal{E}, D(\mathcal{E}))\) is an irreducible transient regular Dirichlet form and \(\mathcal{E}^{\mu}(f,f)\leq \mathcal{E}(f,f)\) holds. For a Bernstein function \(\Phi(\lambda):=1-e^{-\lambda \log2}\), a straightforward computation using the matrix representations of the associated operators shows that \(\mathcal{E}^\Phi(f,f):=\frac{31}{160}(f(1)-f(0))^2+\frac{111}{160}f(1)^2+\frac{129}{320}f(0)^2\) and \(\mathcal{E}^{\mu, \Phi}(f,f)=\frac{15}{32}(f(1)-f(0))^2\). Hence for a function \(f(x)=1_{\{1\}}(x)-1_{\{0\}}(x)\), \(\mathcal{E}^{\mu, \Phi}(f,f)-\mathcal{E}^\Phi(f,f)=1/320>0\) holds.
\end{remark}

We obtain the following representation of the difference between \(\mathcal{E}^{\mu_+, \Phi}\) and \(\mathcal{E}^{\mu, \Phi}\).
\begin{proposition}\label{prop:RR}
It holds that
\begin{eqnarray*}
\mathcal{E}^{\mu_+, \Phi}(f,f) - \mathcal{E}^{\mu, \Phi}(f,f) &=& b\int_E|f|^2\,d\mu_- + \int_0^\infty \langle T_s^\mu f-T_s^{\mu_+} f,f\rangle_m\, d\nu(s)\\
&=& b\int_E|f|^2\,d\mu_- + \int_0^\infty \langle R_\alpha^\mu f-R_{\alpha}^{\mu_+} f,f\rangle_m\, d\widetilde{\nu}(\alpha)
\end{eqnarray*}
for any \(f \in D(\mathcal{E}) \cap C_c.\) Moreover, if \(D(\mathcal{E}^\mu) = D(\mathcal{E}^{\mu_+})\) holds, then
\[\mathcal{E}^{\mu_+, \Phi}(f,f) - \mathcal{E}^{\mu, \Phi}(f,f) = b\int_E|f|^2\,d\mu_- + \int_0^\infty \int_E R_\alpha^\mu f R_\alpha^{\mu_+} f\,d\mu_-\,d\widetilde{\nu}(\alpha)\]
for any \(f \in D(\mathcal{E}) \cap C_c.\)
\end{proposition}
\begin{proof}
By \cite[Theorem 3.2]{OS26}, for \(f\in D(\mathcal{E})\cap C_c(E)\), it holds that
\begin{align*}
\lefteqn{\mathcal{E}^{\mu_+, \Phi}(f,f) - \mathcal{E}^{\mu, \Phi}(f,f)}\\
&= b\,\mathcal{E}^{\mu_+}(f,f) +  \int_0^\infty \langle f-T_s^{\mu_+} f,f\rangle_m\, d\nu(s) - b\,\mathcal{E}^{\mu}(f,f) -  \int_0^\infty \langle f-T_s^{\mu} f,f\rangle_m\, d\nu(s)\\
&=b\int_E|f|^2\,d\mu_- + \int_0^\infty \langle T_s^\mu f-T_s^{\mu_+} f,f\rangle_m\, d\nu(s).
\end{align*}
Furthermore, it holds that
\begin{align}
\nonumber \int_0^\infty \langle T_s^\mu f-T_s^{\mu_+} f,f\rangle_m\, d\nu(s) &=\int_0^\infty \int_0^\infty e^{-\alpha s}\langle T_s^\mu f-T_s^{\mu_+} f,f\rangle_m\, d\widetilde{\nu}(\alpha)\,ds\\
&= \int_0^\infty \langle R_\alpha^\mu f-R_{\alpha}^{\mu_+} f,f\rangle_m\, d\widetilde{\nu}(\alpha). \label{eq:RR_1}
\end{align}

Moreover we assume \(D(\mathcal{E}^\mu) = D(\mathcal{E}^{\mu_+})\) holds. Then, for \(f, g\in D(\mathcal{E})\cap C_c(E),\) we have
\begin{align}
\nonumber \langle R_\alpha^\mu f-R_{\alpha}^{\mu_+} f,g\rangle_m &= \mathcal{E}^{\mu_+}_\alpha( R_\alpha^\mu f-R_\alpha^{\mu_+} f,R_\alpha^{\mu_+}g)\\
\nonumber &=\mathcal{E}^\mu_\alpha(R_\alpha^\mu f, R_\alpha^{\mu_+}g ) + \int_E R_\alpha^\mu f R_\alpha^{\mu_+}g\,d\mu_- - \mathcal{E}^{\mu_+}_\alpha(R_\alpha^{\mu_+} f, R_\alpha^{\mu_+}g )\\
\nonumber &= \langle f, R_\alpha^{\mu_+}g \rangle_m + \int_E R_\alpha^\mu f R_\alpha^{\mu_+}g\,d\mu_- - \langle R_\alpha^{\mu_+} f, g \rangle_m\\
&= \int_E R_\alpha^\mu f R_\alpha^{\mu_+}g\,d\mu_-. \label{eq:RR_2}
\end{align}
Hence the second equality follows from (\ref{eq:RR_1}) and (\ref{eq:RR_2}).
\end{proof}

By using Proposition \ref{prop:RR}, we can define \(\lambda(\mu, \Phi)\) for a complete Bernstein function \(\Phi\) and a closed form \((\mathcal{E}^\mu, D(\mathcal{E})\cap C_c\)) failing the non-negativity by
\[\lambda(\mu, \Phi):=\inf_{0\neq f\in D(\mathcal{E}) \cap C_c(E)} \frac{\mathcal{E}^{\mu_+, \Phi}(f,f)}{ b\int_E|f|^2\,d\mu_- + \int_0^\infty \langle R_\alpha^\mu f-R_{\alpha}^{\mu_+} f,f\rangle_m\, d\widetilde{\nu}(\alpha)}.\]
By definition, \(\lambda(\mu, \Phi)<1\) if and only if \(\mathcal{E}^{\mu, \Phi}\) is supercritical.

\begin{proposition}\label{prop:comp_spectrum}
It holds that \(\lambda(\mu) = \lambda(\mu, \rm{id})\). Moreover, \(\lambda(\mu) \leq \lambda(\mu, \Phi)\) holds for any complete Bernstein function \(\Phi\).
\end{proposition}
\begin{proof}
Since \(\Phi(\lambda)=\lambda\), the drift coefficient \(b=1\) and the Le\'{e}vy measure \(\nu=0\) hold. By Proposition \ref{prop:RR}, \(\lambda(\mu) = \lambda(\mu, \rm{id})\) holds.

Since \(1\leq \lambda(\mu,\Phi)\), the case of \(\lambda(\mu)=1\) is trivial. We assume that \(1<\lambda(\mu)\). For any \(f \in D(\mathcal{E})\cap C_c(E)\), since it holds that
\[\lambda(\mu) \leq \frac{\mathcal{E}^{\mu_+}(f,f)}{\int_E |f|^2\,d\mu_-}.\]
Then we obtain
\[\lambda(\mu)\mathcal{E}^{\mu_+}(f,f) - \lambda(\mu)\mathcal{E}^{\mu}(f,f) = \lambda(\mu)\int_E |f|^2 \,d\mu_- \leq \mathcal{E}^{\mu_+}(f,f)\]
and so
\[\mathcal{E}^{\mu_+}(f,f) \leq \frac{\lambda(\mu)}{\lambda(\mu)-1}\mathcal{E}^{\mu}(f,f).\]
Since \(\Phi\) is a complete Bernstein function, by \cite[Theorem 12.12]{SSV12} and the fact that \(\Phi(c\lambda) \leq c \Phi(\lambda)\) for \(c>1\) and \(\lambda \geq 0\), we have
\[\mathcal{E}^{\mu_+, \Phi}(f,f) \leq \frac{\lambda(\mu)}{\lambda(\mu)-1}\mathcal{E}^{\mu,\Phi}(f,f).\]
Hence we have
\[\lambda(\mu) \leq  \frac{\mathcal{E}^{\mu_+, \Phi}(f,f)}{\mathcal{E}^{\mu_+, \Phi}(f,f) - \mathcal{E}^{\mu, \Phi}(f,f)}\]
and, taking the infimum over all admissible \(f\), we obtain \(\lambda(\mu) \leq \lambda(\mu, \Phi)\). 
\end{proof}

Before considering the relationship between \(\lambda(\mu, \Phi)\) and criticalities of \((\mathcal{E}^{\mu, \Phi}, D(\mathcal{E}^{\mu, \Phi}))\), we give the following proposition concerning well-definedness of the classification of \(\lambda(\mu,\Phi)\) relative to the threshold \(1\).

\begin{proposition}\label{prop:well-def}
For a complete Bernstein function \(\Phi\), Radon measures \(\mu_+, \mu_-\) and \(\xi\), we define
\[\lambda(\mu, \Phi;\xi) := \inf_{0\neq f\in D(\mathcal{E}) \cap C_c(E)} \frac{\mathcal{E}^{\mu_+ + \xi, \Phi}(f,f)}{\mathcal{E}^{\mu_+ + \xi, \Phi}(f,f) - \mathcal{E}^{\mu, \Phi}(f,f)}. \]
Note that \(\lambda(\mu, \Phi;\xi)\) coincides with \(\lambda(\widetilde{\mu}, \Phi)\) as \(\widetilde{\mu}_+:= \mu_+ + \xi\) and \(\widetilde{\mu}_-:= \mu_- + \xi\).
Then \(\lambda(\mu, \Phi;\xi) \leq \lambda(\mu, \Phi)\) holds. In particular, if \(\lambda(\mu, \Phi)=1\) then \(\lambda(\mu, \Phi;\xi)=1\).

Moreover, we assume that the strong Feller property {\rm (SF)} holds for \((\mathcal{E}, D(\mathcal{E}))\) and \(\mu_- + \xi \in \mathcal{K}_\infty(\mathcal{E}^{\mu_+ + \xi})\). Then, \(\lambda(\mu, \Phi)>1\) if and only if \(\lambda(\mu, \Phi;\xi)>1\), and \(\lambda(\mu, \Phi)=1\) if and only if \(\lambda(\mu, \Phi;\xi)=1\).
\end{proposition}

\begin{proof}
Since \(\mathcal{E}^{\mu_+ + \xi, \Phi}(f,f) \geq \mathcal{E}^{\mu_+, \Phi}(f,f)\), it holds that \(\lambda(\mu, \Phi;\xi) \leq \lambda(\mu, \Phi)\).

We assume that \(\mu_- + \xi \in \mathcal{K}_\infty(\mathcal{E}^{\mu_+ + \xi})\) and \(\lambda(\mu, \Phi)>1\). We prove that, there exists \(C>0\) such that, for any \(f\in D(\mathcal{E})\cap C_c\), it holds that
\begin{equation}
\int_E |f|^2 \, d\xi \leq C\,\mathcal{E}^{\mu_+}(f,f) \label{eq:welldef_1}
\end{equation}

by contradiction. Suppose that no such constant exists, then we can take \(f_n \in D(\mathcal{E})\cap C_c\) satisfying \(\int_E |f_n|^2\,d\xi=1\) and \(\mathcal{E}^{\mu_+}(f_n,f_n)\) converges to \(0\). Then \(\{\mathcal{E}^{\mu_+ + \xi}(f_n,f_n)\}_n\) is bounded, so we can take a subsequence \(\{f_{n_k}\}\) converging weakly to some \(u \in D_e(\mathcal{E}^{\mu_+ + \xi})\) in \(\mathcal{E}^{\mu_+ +\xi}\). By \cite[Theorem 1.1 (1)]{KM25}, \((D_e(\mathcal{E}^{\mu_+ + \xi}), \mathcal{E}^{\mu_++\xi})\) is compactly embedded into \(L^2(E;\mu_- + \xi)\). Hence \((D_e(\mathcal{E}^{\mu_+ + \xi}), \mathcal{E}^{\mu_++\xi})\) is compactly embedded into \(L^2(E;\xi)\), and so \(f_{n_k}\) converges to \(u\) in \(L^2(E;\xi)\). Since \(\int_E |f_n|^2\,d\xi=1\), we have \(\int_E |u|^2\,d\xi=1\). However it holds that
\[0 \leq \mathcal{E}^{\mu_+}(u,u) = \mathcal{E}^{\mu_+ + \xi}(u,u)-1 \leq \varliminf_{n _k \to \infty}\left( \mathcal{E}^{\mu_+ + \xi}(f_{n_k},f_{n_k})-1 \right) = \varliminf_{n _k \to \infty}\mathcal{E}^{\mu_+}(f_{n_k},f_{n_k})= 0.\]
By the transience of \(\mathcal{E}^{\mu_+}\), \(u=0\) holds. This is a contradiction to  \(\int_E |u|^2\,d\xi=1\).

By (\ref{eq:welldef_1}), we have \(\mathcal{E}^{\mu_+ + \xi}(f,f) \leq (1+C) \mathcal{E}^{\mu_+}(f,f)\). Since \(\Phi\) is complete and concave,  \(\mathcal{E}^{\mu_+ + \xi, \Phi}(f,f) \leq C_1 \mathcal{E}^{\mu_+, \Phi}(f,f)\) holds for some \(C_1>1\) independent of \(f\). Hence we have
\begin{eqnarray*}
\frac{\mathcal{E}^{\mu_+ + \xi, \Phi}(f,f)}{\mathcal{E}^{\mu_+ + \xi, \Phi}(f,f) - \mathcal{E}^{\mu, \Phi}(f,f)} & \geq & \frac{C_1 \mathcal{E}^{\mu_+, \Phi}(f,f)}{C_1 \mathcal{E}^{\mu_+, \Phi}(f,f) - \mathcal{E}^{\mu, \Phi}(f,f)}\\
& \geq & \frac{C_1 \mathcal{E}^{\mu_+, \Phi}(f,f)}{C_1  \mathcal{E}^{\mu_+, \Phi}(f,f) - \frac{\lambda(\mu,\Phi)-1}{\lambda(\mu,\Phi)} \mathcal{E}^{\mu_+, \Phi}(f,f)}\\
&=& \frac{C_1}{C_1  - \frac{\lambda(\mu,\Phi)-1}{\lambda(\mu,\Phi)}}.
\end{eqnarray*}
Since \(\lambda(\mu, \Phi)>1\), it holds that \(\lambda(\mu, \Phi;\xi)>1\).
\end{proof}

\begin{remark}
We give an example satisfying \(\lambda(\mu, \Phi)>1\) but \(\lambda(\mu, \Phi;\xi)=1\).
Let \(\mathcal{E}(f,f):=\int_{\mathbb{R}}|\nabla f|^2\,dx + |f(0)|^2\) and \(D(\mathcal{E}):=H^1(\mathbb{R})\). Then \((\mathcal{E}, D(\mathcal{E}))\) is an irreducible transient regular Dirichlet form on \(L^2(\mathbb{R})\). We set \(\mu_+:=0, \mu_-:= \frac{1}{2}\delta_0\) and \(\Phi={\rm id}.\) Then we have \(\mathcal{E}^\mu(f,f) = \int_{\mathbb{R}}|\nabla f|^2\,dx + \frac{1}{2}|f(0)|^2 \) and 
\[\lambda(\mu)=\inf_{f\in D(\mathcal{E}) \cap C_c(E)} \frac{\int_{\mathbb{R}}|\nabla f|^2\,dx + |f(0)|^2}{\frac{1}{2}|f(0)|^2 } \geq 2. \]
However, for \(d\xi:=dx\), we have \(\mathcal{E}^\xi(f,f)=\int_{\mathbb{R}}|\nabla f|^2\,dx + |f(0)|^2 + \int_{\mathbb{R}}|f|^2\,dx\). We fix a non-zero function \(\varphi \in C_c^\infty\) whose support is in \((-1,1)\) and we set \(f_n(x):=\varphi((x-2n)/n)\) for each \(n \in \mathbb{N}\). Then we have
\[\mathcal{E}^{\xi}(f_n,f_n) = \frac{1}{n} \int |\nabla \varphi|^2\,dx + 0 + n\int|\varphi|^2\,dx\]
and
\[\int |f_n|^2\,d(\mu_- + \xi) =n\int|\varphi|^2\,dx. \]
Hence we have
\[1\leq \lambda(\mu, \Phi;\xi) \leq \frac{\mathcal{E}^{\xi}(f_n,f_n)}{\int |f_n|^2\,d(\mu_- + \xi)} = 1+ \frac{1}{n^2}\frac{\int|\nabla \varphi|^2\,dx}{\int|\varphi|^2\,dx} \to 1.\]
\end{remark}

The above example shows that, even for Schr\"{o}dinger forms, the distinction between \(\lambda=1\) and \(\lambda>1\) is not necessarily preserved when the same positive measure \(\xi\) is added to both \(\mu_+\) and \(\mu_-\). However, within Takeda's standard framework for the criticality of Schr\"{o}dinger forms, this distinction is independent of the choice of decomposition, as long as the corresponding assumptions on the Green-tightness are satisfied. Moreover, the above example shows that \(\lambda=1\) does not imply the criticality in general.

We also remark that  \(\lambda(\mu, \Phi;\xi)>1\) if and only if \(\lambda(\mu, \Phi)>1\) for \(\xi \in \mathcal{K}_\infty(\mathcal{E}^{\mu_+})\) by a similar proof of Proposition \ref{prop:well-def} with the Stollmann-Voigt inequality \cite{SV96}.

We establish a criterion for subcriticality and criticality by using \(\lambda(\mu, \Phi)\). From the following theorem, we obtain a sufficient condition for subcriticality
\begin{theorem}\label{thm:subcri}
For a signed smooth Radon measure \(\mu\) and a complete Bernstein function \(\Phi\), if \(\lambda(\mu, \Phi)>1\) then \((\mathcal{E}^{\mu, \Phi}, \mathcal{D}(\mathcal{E}^{\mu, \Phi}))\) is subcritical.
\end{theorem}
\begin{proof}
As in the proof of Proposition \ref{prop:comp_spectrum}, we have
\begin{equation}
\mathcal{E}^{\mu_+, \Phi}(f,f) \leq \frac{\lambda(\mu,\Phi)}{\lambda(\mu,\Phi)-1}\mathcal{E}^{\mu,\Phi}(f,f) \label{eq:subcri_1}
\end{equation}
for any \(f\in D(\mathcal{E})\cap C_c\). Since we assumed \((\mathcal{E}^{\mu_+}, D(\mathcal{E}^{\mu_+}))\) is transient, by \cite[Theorem 3.2]{O02}, \((\mathcal{E}^{\mu_+,\Phi}, D(\mathcal{E}^{\mu_+,\Phi}))\) is also a transient Dirichlet form. By \cite[Theorem 2.1.5 (i)]{CF12},  there exists a strictly positive bounded function \(g\in L^1(E;m)\) satisfying
\[\int_E |f|g \,dm \leq \sqrt{\mathcal{E}^{\mu_+,\Phi}(f,f)}\]
for any \(f\in D(\mathcal{E}^{\mu_+,\Phi})\). Combining this with (\ref{eq:subcri_1}), we have 
\begin{equation}
\sqrt{\frac{\lambda(\mu,\Phi)-1}{\lambda(\mu,\Phi)}} \int_E |f|g \,dm \leq  \sqrt{\mathcal{E}^{\mu,\Phi}(f,f)} \label{eq:subcri_2}
\end{equation}
for any \(f\in D(\mathcal{E})\cap C_c\). The inequality (\ref{eq:subcri_2}) holds for \(f\in D_e(\mathcal{E}^{\mu,\Phi})\), and so \((\mathcal{E}^{\mu, \Phi}, \mathcal{D}(\mathcal{E}^{\mu, \Phi}))\) is subcritical.
\end{proof}

Before establishing a criterion for criticality, we state the following lemmas.
\begin{lemma}\label{lem:Res_conti}
For a Radon measure \(\mu_+\), a non-zero Bernstein function \(\Phi\) and each \(\alpha>0\), a linear operator \(R_\alpha^{\mu_+} : (D_e(\mathcal{E}^{\mu_+, \Phi}), \mathcal{E}^{\mu_+, \Phi}) \to (D_e(\mathcal{E}^{\mu_+}), \mathcal{E}^{\mu_+})\) is continuous.
\end{lemma}
\begin{proof}
First, we consider the case of \(b>0\). In this case, by \cite[Theorem 2.1 (iv)]{O02}, \((D_e(\mathcal{E}^{\mu_+, \Phi}), \mathcal{E}^{\mu_+, \Phi})\) is continuously embedded into \((D_e(\mathcal{E}^{\mu_+}), \mathcal{E}^{\mu_+})\). By the spectral decomposition, for each \(\alpha>0\), a linear operator \(R_\alpha^{\mu_+}:(D_e(\mathcal{E}^{\mu_+}), \mathcal{E}^{\mu_+}) \to (D_e(\mathcal{E}^{\mu_+}), \mathcal{E}^{\mu_+})\) is continuous, and so is  \(R_\alpha^{\mu_+} : (D_e(\mathcal{E}^{\mu_+, \Phi}), \mathcal{E}^{\mu_+, \Phi}) \to (D_e(\mathcal{E}^{\mu_+}), \mathcal{E}^{\mu_+})\).

Next we consider the case of \(b=0\). Without loss of generality, we may assume that \(\nu \not =0\). For any \(f \in D(\mathcal{E}^{\mu_+, \Phi}) \subset L^2(E;m)\), we have
\begin{eqnarray}
\nonumber \alpha \mathcal{E}^{\mu_+}(R_\alpha^{\mu_+} f, R_\alpha^{\mu_+} f) &=& \langle f- \alpha R_\alpha^{\mu_+} f, \alpha R_\alpha^{\mu_+} f \rangle_m \\
\nonumber &\leq &  \langle f- \alpha R_\alpha^{\mu_+} f, \alpha R_\alpha^{\mu_+} f \rangle_m + \|f- \alpha R_\alpha^{\mu_+} f\|^2_{L^2(E;m)}\\
\nonumber &=&  \langle f- \alpha R_\alpha^{\mu_+} f, f \rangle_m \\
&=& \int_0^\infty \alpha e^{-\alpha s} \langle f- T_s^{\mu_+} f, f \rangle_m \,ds.\label{eq:lem_Res_conti_1}
\end{eqnarray}
Since \(\langle f- T_s^{\mu_+} f, f \rangle_m\) is increasing and \(\frac{1}{s}\langle f- T_s^{\mu_+} f, f \rangle_m\) is decreasing, for any \(t>0\), we have
\begin{eqnarray}
\nonumber \int_0^\infty \alpha e^{-\alpha s} \langle f- T_s^{\mu_+} f, f \rangle_m \,ds &\leq & \langle f-  T_t^{\mu_+} f, f \rangle_m \left( \int_0^t \alpha e^{-\alpha s}  \,ds + \int_t^\infty \frac{s}{t}\alpha e^{-\alpha s}\,ds \right)\\
&=& \langle f-  T_t^{\mu_+} f, f \rangle_m  \left( 1+\frac{e^{-\alpha t}}{\alpha t} \right).\label{eq:lem_Res_conti_2}
\end{eqnarray}

By \cite[Theorem 2.1(i), (v)]{O02}, we have
\begin{eqnarray}
\nonumber \mathcal{E}^{\mu_+, \Phi}(f,f) &=& \int_0^\infty \langle f-  T_s^{\mu_+} f, f \rangle_m \,d\nu(s)\\
\nonumber &\geq &  \int_t^\infty \langle f-  T_s^{\mu_+} f, f \rangle_m \,d\nu(s)\\
&\geq & \langle f-  T_t^{\mu_+} f, f \rangle_m \nu([t,\infty)).\label{eq:lem_Res_conti_3}
\end{eqnarray}
Since \(\nu\) is a L\'evy measure on \((0,\infty)\), \( \nu([t,\infty))\) is finite. Then, by (\ref{eq:lem_Res_conti_1}), (\ref{eq:lem_Res_conti_2}) and (\ref{eq:lem_Res_conti_3}), for large \(t>0\), it holds that
\begin{equation}
\mathcal{E}^{\mu_+}(R_\alpha^{\mu_+} f, R_\alpha^{\mu_+} f)  \leq \frac{1}{\alpha}\left( 1+\frac{e^{-\alpha t}}{\alpha t} \right) \frac{1}{\nu([t,\infty))} \mathcal{E}^{\mu_+, \Phi}(f,f). \label{eq:lem_Res_conti_4}
\end{equation}

For \(f \in D_e(\mathcal{E}^{\mu_+, \Phi}),\) by taking an approximating sequence in \(D(\mathcal{E}^{\mu_+, \Phi})\), (\ref{eq:lem_Res_conti_4}) holds on \(D_e(\mathcal{E}^{\mu_+, \Phi})\). Hence  \(R_\alpha^{\mu_+} : (D_e(\mathcal{E}^{\mu_+, \Phi}), \mathcal{E}^{\mu_+, \Phi}) \to (D_e(\mathcal{E}^{\mu_+}), \mathcal{E}^{\mu_+})\) is continuous.
\end{proof}

\begin{remark}\label{rem_lem1}
Unlike the other results, Lemma \ref{lem:Res_conti} only requires that \(\Phi\) be a Bernstein function, with no additional assumptions. In particular, for any Bernstein function \(\Phi\), there exists \(C>0\) depending only on \(\Phi\) such that, for any \(f\in D_e(\mathcal{E}^{\mu_+, \Phi}),\) it holds that
\[\mathcal{E}^{\mu_+}(R_\alpha^{\mu_+} f, R_\alpha^{\mu_+} f)  \leq \frac{C}{\alpha}\left( 1+\frac{1}{\alpha} \right) \mathcal{E}^{\mu_+, \Phi}(f,f).\]
\end{remark}

For a signed smooth Radon measure \(\mu \in \mathcal{K}_{loc}(\mathcal{E})- \mathcal{K}_{\infty}(\mathcal{E}^{\mu_+})\), we define \(\theta_\alpha\) by
\begin{eqnarray*}
\theta_\alpha := \sup_{0\neq f \in D(\mathcal{E}^{\mu_+})} \frac{\sqrt{\int_E |f|^2\, d\mu_-}}{\sqrt{\mathcal{E}_\alpha^{\mu_+}(f,f)}} < 1.
\end{eqnarray*}
The above inequality follows from the compact embedding of \((D_e(\mathcal{E}^{\mu_+}), \mathcal{E}^{\mu_+})\) into \(L^2(E;\mu_-)\) (See \cite[Theorem 1.1 (1)]{KM25} for example) and the non-negativity of \(\mathcal{E}^\mu\).

\begin{lemma}\label{lem:Res_conti2}
Suppose the strong Feller property {\rm (SF)} holds for \((\mathcal{E}^{\mu_+}, D(\mathcal{E}^{\mu_+}))\). Then, for a signed smooth Radon measure \(\mu \in \mathcal{K}_{loc}(\mathcal{E})-\mathcal{K}_{\infty}(\mathcal{E}^{\mu_+})\), a Bernstein function \(\Phi\), each \(\alpha>0\) and any \(f\in D_e(\mathcal{E}^{\mu_+, \Phi})\), it holds that
\begin{equation}
\|R_\alpha^\mu f\|_{L^2(\mu_-)} \leq \frac{1}{1-\theta_\alpha} \| R_\alpha^{\mu_+} f \|_{L^2(\mu_-)}. \label{eq:critical_1_2}
\end{equation}
and
\begin{equation}
\mathcal{E}^{\mu_+}(R_\alpha^{\mu}f,R_\alpha^{\mu}f) \leq  \left(\frac{1}{1-\theta_\alpha}\right)^2 \mathcal{E}^{\mu_+}(R_\alpha^{\mu_+} f,R_\alpha^{\mu_+} f).\label{eq:critical_1_a}
\end{equation}
In particular, an linear operator \(R_\alpha^{\mu} : (D_e(\mathcal{E}^{\mu_+, \Phi}), \mathcal{E}^{\mu_+, \Phi}) \to (D_e(\mathcal{E}^{\mu_+}), \mathcal{E}^{\mu_+})\) is continuous.
\end{lemma}
\begin{proof}
The identity \(D(\mathcal{E}^{\mu_+})=D(\mathcal{E}^{\mu})\) follows from the Stollmann-Voigt inequality \cite{SV96} (see also \cite[Theorem 2.8]{OS26}), and so, it holds that \(R_\alpha^\mu f \in D(\mathcal{E}^{\mu_+})=D(\mathcal{E}^{\mu})\) for any \(f\in D(\mathcal{E}^{\mu_+, \Phi}) \subset L^2(E;m)\). 

By the Cauchy--Schwarz inequality and the definition of \(\theta_\alpha\), we have
\begin{eqnarray*}
\lefteqn{\mathcal{E}_\alpha^{\mu_+}(R_\alpha^{\mu}f-R_\alpha^{\mu_+} f, R_\alpha^{\mu}f-R_\alpha^{\mu_+} f)}\\
&=&  \mathcal{E}_\alpha^{\mu_+}(R_\alpha^{\mu}f, R_\alpha^{\mu}f-R_\alpha^{\mu_+} f) - \mathcal{E}_\alpha^{\mu_+}(R_\alpha^{\mu_+}f, R_\alpha^{\mu}f-R_\alpha^{\mu_+} f)\\
&=&  \mathcal{E}_\alpha^{\mu}(R_\alpha^{\mu}f, R_\alpha^{\mu}f-R_\alpha^{\mu_+} f) +\int_E R_\alpha^{\mu}f \left(R_\alpha^{\mu}f-R_\alpha^{\mu_+} f\right)\,d\mu_- - \langle f, R_\alpha^{\mu}f-R_\alpha^{\mu_+} f\rangle_m\\
&=& \int_E R_\alpha^{\mu}f \left(R_\alpha^{\mu}f-R_\alpha^{\mu_+} f\right)\,d\mu_-\\
&\leq & \| R_\alpha^\mu f \|_{L^2(\mu_-)}\, \| R_\alpha^{\mu}f-R_\alpha^{\mu_+}f\|_{L^2(\mu_-)}\\
&\leq & \| R_\alpha^\mu f \|_{L^2(\mu_-)}\, \theta_\alpha \, \sqrt{\mathcal{E}_\alpha^{\mu_+}(R_\alpha^{\mu}f-R_\alpha^{\mu_+} f,R_\alpha^{\mu}f-R_\alpha^{\mu_+} f)}
\end{eqnarray*}
Hence we have
\begin{equation}
\sqrt{\mathcal{E}_\alpha^{\mu_+}(R_\alpha^{\mu}f-R_\alpha^{\mu_+} f,R_\alpha^{\mu}f-R_\alpha^{\mu_+} f)} \leq  \theta_\alpha \, \| R_\alpha^\mu f \|_{L^2(\mu_-)} \label{eq:critical_1_1}
\end{equation}
Combining (\ref{eq:critical_1_1}) with the inequality
\begin{eqnarray*}
\|R_\alpha^\mu f\|_{L^2(\mu_-)} - \|R_\alpha^{\mu_+} f\|_{L^2(\mu_-)} \leq \|R_\alpha^\mu f - R_\alpha^{\mu_+} f\|_{L^2(\mu_-)} \leq \sqrt{\mathcal{E}_\alpha^{\mu_+}(R_\alpha^{\mu}f-R_\alpha^{\mu_+} f,R_\alpha^{\mu}f-R_\alpha^{\mu_+} f)},
\end{eqnarray*}
we obtain
\begin{equation*}
\|R_\alpha^\mu f\|_{L^2(\mu_-)} \leq \frac{1}{1-\theta_\alpha} \| R_\alpha^{\mu_+} f \|_{L^2(\mu_-)}.
\end{equation*}

By (\ref{eq:critical_1_1}), (\ref{eq:critical_1_2}) and the non-negativity of \(\mathcal{E}^\mu\), we have
\begin{eqnarray}
\nonumber \sqrt{\mathcal{E}^{\mu_+}(R_\alpha^{\mu}f,R_\alpha^{\mu}f)} &\leq &\sqrt{\mathcal{E}^{\mu_+}(R_\alpha^{\mu}f-R_\alpha^{\mu_+} f,R_\alpha^{\mu}f-R_\alpha^{\mu_+} f)} +\sqrt{\mathcal{E}^{\mu_+}(R_\alpha^{\mu_+} f,R_\alpha^{\mu_+} f)} \\
\nonumber &\leq & \theta_\alpha \, \| R_\alpha^\mu f \|_{L^2(\mu_-)} + \sqrt{\mathcal{E}^{\mu_+}(R_\alpha^{\mu_+} f,R_\alpha^{\mu_+} f)}\\
\nonumber &\leq &  \frac{\theta_\alpha}{1-\theta_\alpha} \| R_\alpha^{\mu_+} f \|_{L^2(\mu_-)} + \sqrt{\mathcal{E}^{\mu_+}(R_\alpha^{\mu_+} f,R_\alpha^{\mu_+} f)}\\
&\leq &  \frac{1}{1-\theta_\alpha}  \sqrt{\mathcal{E}^{\mu_+}(R_\alpha^{\mu_+} f,R_\alpha^{\mu_+} f)}. \label{eq:critical_1_3}
\end{eqnarray}
For \(f \in D_e(\mathcal{E}^{\mu_+, \Phi}),\) by taking an approximating sequence in \(D(\mathcal{E}^{\mu_+, \Phi})\), (\ref{eq:critical_1_3}) holds on \(D_e(\mathcal{E}^{\mu_+, \Phi})\).

Combining these with Lemma \ref{lem:Res_conti}, an linear operator \(R_\alpha^{\mu} : (D_e(\mathcal{E}^{\mu_+, \Phi}), \mathcal{E}^{\mu_+, \Phi}) \to (D_e(\mathcal{E}^{\mu_+}), \mathcal{E}^{\mu_+})\) is continuous.
\end{proof}

\begin{remark}
Although we need to assume that {\rm (SF)}, \(\mu_+ \in \mathcal{K}_{loc}(\mathcal{E})\) and \(\mu_-\in \mathcal{K}_{\infty}(\mathcal{E}^{\mu_+})\), Lemma \ref{lem:Res_conti2} also only requires that \(\Phi\) be a Bernstein function, with no additional assumptions. Combining Lemma \ref{lem:Res_conti} with Lemma \ref{lem:Res_conti2}, for any Bernstein function \(\Phi\), there exists \(C>0\) depending only on \(\Phi\) such that, for any \(f\in D_e(\mathcal{E}^{\mu_+, \Phi}),\)
\[\mathcal{E}^{\mu_+}(R_\alpha^{\mu} f, R_\alpha^{\mu} f)  \leq \frac{C}{\alpha}\left( 1+\frac{1}{\alpha} \right) \left(\frac{1}{1-\theta_\alpha}\right)^2 \mathcal{E}^{\mu_+, \Phi}(f,f).\]
\end{remark}

By using the Stollmann-Voigt inequality, we obtain the following estimate.
\begin{lemma}\label{lem:Ehrling's inequality}
Suppose that \(\mu_- \in \mathcal{K}_\infty(\mathcal{E}^{\mu_+})\). Then, for any \(\varepsilon>0\), there exists \(C_\varepsilon >0\) such that, for any \(f \in D(\mathcal{E}^{\mu_+})\),
\begin{equation}
\int_E |f|^2\,d\mu_- \leq \varepsilon \mathcal{E}^{\mu_+}(f,f) + C_\varepsilon \mathcal{E}^{\mu_+}(R_1^{\mu_+}f,R_1^{\mu_+}f). \label{eq:Ehrling's inequality}
\end{equation}
\end{lemma}
\begin{proof}
We note that the Green-tight class \(\mathcal{K}_\infty(\mathcal{E}^{\mu_+})\) is included in the Kato class \(\mathcal{K}(\mathcal{E}^{\mu_+})\) by the definition.

We fix \(\varepsilon >0\). Since \(\mu_- \in \mathcal{K}(\mathcal{E}^{\mu_+})\), we can choose large \(\beta>1\) satisfying \(\|R_\beta^{\mu_+}\mu_-\|_\infty \leq \varepsilon/2\).
By the Stollmann-Voigt inequality for \((D(\mathcal{E}^{\mu_+}), \mathcal{E}^{\mu_+}_\beta)\) and the contraction property of the resolvent, for any \(f \in D(\mathcal{E}^{\mu_+})\), we have
\begin{eqnarray}
\nonumber \int_E |f-\beta R_\beta^{\mu_+}f|^2 \,d\mu_- &\leq & \|R_\beta^{\mu_+}\mu_-\|_\infty \mathcal{E}_\beta^{\mu_+}(f-\beta R_\beta^{\mu_+}f,f-\beta R_\beta^{\mu_+}f)\\
 &\leq & \frac{\varepsilon}{2}\, \mathcal{E}^{\mu_+}(f, f) \label{eq:SVkai_1}
\end{eqnarray}
By the Stollmann-Voigt inequality for \((D(\mathcal{E}^{\mu_+}), \mathcal{E}^{\mu_+})\) 
\begin{eqnarray}
\nonumber \int_E |\beta R_\beta^{\mu_+}f|^2 \,d\mu_- &\leq & \|R^{\mu_+}\mu_-\|_\infty \mathcal{E}^{\mu_+}(\beta R_\beta^{\mu_+}f,\beta R_\beta^{\mu_+}f)\\
\nonumber &=& \|R^{\mu_+}\mu_-\|_\infty \beta^2 \mathcal{E}^{\mu_+}( R_\beta^{\mu_+}f, R_\beta^{\mu_+}f)\\
&\leq &\|R^{\mu_+}\mu_-\|_\infty \beta^2 \mathcal{E}^{\mu_+}( R_1^{\mu_+}f, R_1^{\mu_+}f).  \label{eq:SVkai_2}
\end{eqnarray}
Here we used the inequality \(\mathcal{E}^{\mu_+}(R_\beta^{\mu_+}f, R_\beta^{\mu_+}f) \leq \mathcal{E}^{\mu_+}(R_1^{\mu_+}f, R_1^{\mu_+}f)\) obtained by the spectral decomposition.

Hence, by setting \(C_\varepsilon := 2 \|R^{\mu_+}\mu_-\|_\infty \beta^2\), we obtain the inequality (\ref{eq:Ehrling's inequality}).
\end{proof}

We introduce the following uniform convergence condition {\rm (UC)}.
\begin{equation}
\varlimsup_{\delta \searrow 0}\sup \left\{\int_0^\delta \int_E  R_\alpha^\mu f\, R_\alpha^{\mu_+} f\,d\mu_-\,d\widetilde{\nu}(\alpha) : f \in D(\mathcal{E}^{\mu_+, \Phi}) \text{\ with\ } \mathcal{E}^{\mu_+, \Phi}(f,f)\leq 1 \right\} =0. \tag{UC}
\end{equation}

\begin{theorem}\label{thm:cri}
Suppose that the strong Feller property {\rm (SF)} holds for \((\mathcal{E}^{\mu_+}, D(\mathcal{E}^{\mu_+}))\). For a signed smooth Radon measure \(\mu \in \mathcal{K}_{loc}(\mathcal{E})- \mathcal{K}_{\infty}(\mathcal{E}^{\mu_+})\) and a complete Bernstein function \(\Phi\) satisfying {\rm (UC)}, if \(\lambda(\mu, \Phi)=1\) then \((\mathcal{E}^{\mu, \Phi}, \mathcal{D}(\mathcal{E}^{\mu, \Phi}))\) is critical.
\end{theorem}
\begin{proof}
Without loss of generality, we may assume that \(\Phi \neq 0.\) Since \(\lambda(\mu, \Phi)=1\), there exists a sequence \(\{f_n\}_n \subset D(\mathcal{E})\cap C_c\) satisfying \(\mathcal{E}^{\mu_+, \Phi}(f_n,f_n)-\mathcal{E}^{\mu, \Phi}(f_n,f_n) = 1 \) and \(\mathcal{E}^{\mu_+, \Phi}(f_n,f_n)\) converges to \(1\) as \(n \to \infty\). Since \(\{f_n\}_n\) is a bounded sequence in \((D_e(\mathcal{E}^{\mu_+, \Phi}), \mathcal{E}^{\mu_+, \Phi})\), there exists a subsequence \(\{n_k\}\) such that \(f_{n_k}\) converges weakly to some function \(v \in D_e(\mathcal{E}^{\mu_+, \Phi}) \) in  \(\mathcal{E}^{\mu_+, \Phi}\).

By Lemma \ref{lem:Res_conti}, \(R_\alpha^{\mu_+} : (D_e(\mathcal{E}^{\mu_+, \Phi}), \mathcal{E}^{\mu_+, \Phi}) \to (D_e(\mathcal{E}^{\mu_+}), \mathcal{E}^{\mu_+})\) is continuous. By \cite[Theorem 1.1 (1)]{KM25}, \((D_e(\mathcal{E}^{\mu_+}), \mathcal{E}^{\mu_+})\) is compactly embedded into \(L^2(E;\mu_-)\). We remark that, although a Green-tight measure in \cite[Theorem 1.1 (1)]{KM25} is a one in the sense of Chen \cite{C02}, and a Green-tight measure in Definition \ref{defKato} is one in the sense of Zhao \cite{Z92}, both classes coincide under the condition (SF) by \cite[Lemma 4.1]{KK17} (or \cite[Theorem 1.2]{KM25}). Hence, the linear operator \({\rm id} \circ R_\alpha^{\mu_+}: (D_e(\mathcal{E}^{\mu_+, \Phi}),\mathcal{E}^{\mu_+, \Phi}) \to L^2(E;\mu_-)\) is a compact operator, and so  \(R_\alpha^{\mu_+} f_{n_k}\) converges strongly to \(R_\alpha^{\mu_+} v \in D_e(\mathcal{E}^{\mu_+})\) in  \(L^2(E;\mu_-)\) for each \(\alpha >0\). 

Similarly, by Lemma \ref{lem:Res_conti2}, the linear operator \({\rm id} \circ R_\alpha^{\mu}: (D_e(\mathcal{E}^{\mu_+, \Phi}),\mathcal{E}^{\mu_+, \Phi}) \to L^2(E;\mu_-)\) is a compact operator, and so  \(R_\alpha^{\mu} f_{n_k}\) converges strongly to \(R_\alpha^{\mu} v \in D_e(\mathcal{E}^{\mu_+})\) in  \(L^2(E;\mu_-)\) for each \(\alpha >0\).

By Proposition \ref{prop:RR}, we have
\begin{equation}
1 = \mathcal{E}^{\mu_+, \Phi}(f_{n_k},f_{n_k}) - \mathcal{E}^{\mu, \Phi}(f_{n_k},f_{n_k})= b\int_E|f_{n_k}|^2\,d\mu_- + \int_0^\infty \int_E R_\alpha^\mu f_{n_k}\, R_\alpha^{\mu_+} f_{n_k}\,d\mu_-\,d\widetilde{\nu}(\alpha). \label{eq:critical_1_4}
\end{equation}

By the assumption {\rm (UC)}, for a small \(\varepsilon>0\), we take \(\delta >0\) satisfying
\[\int_0^\delta \int_E  R_\alpha^\mu f\, R_\alpha^{\mu_+} f\,d\mu_-\,d\widetilde{\nu}(\alpha) \leq \varepsilon \]
for any \(f \in D_e(\mathcal{E}^{\mu_+, \Phi})\) satisfying \(\mathcal{E}^{\mu_+, \Phi}(f,f) \leq 1.\) Since \(\{f_{n_k}\}\) is bounded with respect to \(\mathcal{E}^{\mu_+, \Phi}\), there exists \(C>0\) such that
\begin{equation}
\int_0^\delta \int_E  R_\alpha^\mu f_{n_k}\, R_\alpha^{\mu_+}f_{n_k}\,d\mu_-\,d\widetilde{\nu}(\alpha) \leq C \varepsilon. \label{eq:critical_1_5}
\end{equation}

We now prove the convergence of \(\mathcal{E}^{\mu_+, \Phi}(f_{n_k},f_{n_k}) - \mathcal{E}^{\mu, \Phi}(f_{n_k},f_{n_k})\) by considering two cases, \(b>0\) and \(b=0\).

\medskip
\noindent\textbf{Case 1: \(b>0\)}

We consider the case \(b>0\). By (\ref{eq:critical_1_2}) and the non-negativity of \(\mathcal{E}^\mu\), we have
\begin{eqnarray*}
\int_E  R_\alpha^\mu f_{n_k}\, R_\alpha^{\mu_+}f_{n_k}\,d\mu_- &\leq &\frac{1}{1-\theta_\alpha} \| R_\alpha^{\mu_+} f_{n_k} \|_{L^2(\mu_-)}^2\\
 &\leq &\frac{1}{1-\theta_\alpha} \mathcal{E}^{\mu_+}(R_\alpha^{\mu_+} f_{n_k},R_\alpha^{\mu_+} f_{n_k}).
\end{eqnarray*}

Since \(\lim_{n\to \infty}\mathcal{E}^{\mu_+, \Phi}(f_n,f_n)=1\) and \(\theta_\alpha\) is non-increasing, it holds that, for any \(\delta \leq \alpha\),
\begin{equation}
\int_E  R_\alpha^\mu f_{n_k}\, R_\alpha^{\mu_+}f_{n_k}\,d\mu_- \leq \frac{1}{1-\theta_\delta} \frac{1}{\alpha^2} \frac{1}{b} \mathcal{E}^{\mu_+, \Phi}(f_{n_k}, f_{n_k}) \leq \frac{C_\delta}{\alpha^2}  \label{eq:critical_1_6}
\end{equation}
for some constant \(C_\delta>0\) depending only on \(\delta\).
By (\ref{eq:Levytilde}), it holds that
\[\int_\delta^\infty \frac{C_\delta}{\alpha^2} \,d\widetilde{\nu}(\alpha) \leq C_\delta \frac{\delta +1}{\delta} \int_\delta^\infty  \frac{1}{\alpha (\alpha+1)}  \,d\widetilde{\nu}(\alpha) < \infty, \]
and so the dominated convergence theorem yields that
\begin{equation}
\lim_{n_k \to \infty}\int_\delta ^\infty \int_E  R_\alpha^\mu f_{n_k}\, R_\alpha^{\mu_+}f_{n_k}\,d\mu_-\,d\widetilde{\nu}(\alpha) = 
\int_\delta ^\infty \int_E  R_\alpha^\mu v\, R_\alpha^{\mu_+}v\,d\mu_-\,d\widetilde{\nu}(\alpha). \label{eq:critical_1_7}
\end{equation}
By the inequality \(b\,\mathcal{E}^{\mu_+}(f_{n_k},f_{n_k}) \leq \mathcal{E}^{\mu_+, \Phi}(f_{n_k},f_{n_k})\) and the compact embedding \(D_e(\mathcal{E}^{\mu_+})\) into \(L^2(E;\mu_-)\), \(f_{n_k}\) converges to \(v\) in \(L^2(E;\mu_-)\). In particular,  \(b\int_E|f_{n_k}|^2\,d\mu_- \) converges to \(b\int_E|v|^2\,d\mu_- \). 

Combining these with (\ref{eq:critical_1_5}), \(\mathcal{E}^{\mu_+, \Phi}(f_{n_k},f_{n_k}) - \mathcal{E}^{\mu, \Phi}(f_{n_k},f_{n_k})\) converges to \(\mathcal{E}^{\mu_+, \Phi}(v,v) - \mathcal{E}^{\mu, \Phi}(v,v)\) for the case of \(b>0\).

\medskip
\noindent\textbf{Case 2: \(b=0\)} 

We consider the case of \(b=0\). For any \(\varepsilon>0\), by applying Lemma \ref{lem:Ehrling's inequality} to the functions \(R_\alpha^{\mu_+}f_{n_k}\), we have
\begin{equation}
\int_E |R_\alpha^{\mu_+}f_{n_k}|^2\,d\mu_- \leq \varepsilon \mathcal{E}^{\mu_+}(R_\alpha^{\mu_+}f_{n_k},R_\alpha^{\mu_+}f_{n_k}) + C_\varepsilon \mathcal{E}^{\mu_+}(R_1^{\mu_+}R_\alpha^{\mu_+}f_{n_k}, R_1^{\mu_+}R_\alpha^{\mu_+}f_{n_k}). \label{eq:E_1}
\end{equation}
We take \(C_1:=\sup_{n_k}\mathcal{E}^{\mu_+, \Phi}(f_{n_k},f_{n_k}) <\infty,\) and we fix large \(M>1\) satisfying \(\int_M^\infty \alpha^{-2} \,d\widetilde{\nu}(\alpha) \leq \varepsilon/C_\varepsilon .\) Then, by the Cauchy–Schwarz inequality, Lemma \ref{lem:Res_conti2}, and (\ref{eq:E_1}), we have
\begin{eqnarray}
\nonumber \lefteqn{\int_M ^\infty \int_E R_\alpha^{\mu_+}f_{n_k}R_\alpha^{\mu}f_{n_k}\,d\mu_- \,d\widetilde{\nu}(\alpha)}\\
\nonumber  &\leq &\int_M ^\infty \frac{1}{1-\theta_\alpha} \int_E |R_\alpha^{\mu_+}f_{n_k}|^2\,d\mu_- \,d\widetilde{\nu}(\alpha)\\
\nonumber &\leq &\frac{\varepsilon}{1-\theta_1}  \int_M ^\infty  \mathcal{E}^{\mu_+}(R_\alpha^{\mu_+}f_{n_k},R_\alpha^{\mu_+}f_{n_k})\,d\widetilde{\nu}(\alpha)\\
&&\quad + \frac{C_\varepsilon }{1-\theta_1}  \int_M ^\infty  \mathcal{E}^{\mu_+}(R_1^{\mu_+}R_\alpha^{\mu_+}f_{n_k}, R_1^{\mu_+}R_\alpha^{\mu_+}f_{n_k})\,d\widetilde{\nu}(\alpha)\label{eq:cri_b_1}
\end{eqnarray}

By (\ref{eq:lem_Res_conti_1}), we have
\begin{eqnarray}
\nonumber \int_0^\infty \mathcal{E}^{\mu_+}(R_\alpha^{\mu_+} f_{n_k}, R_\alpha^{\mu_+} f_{n_k}) \,d\widetilde{\nu}(\alpha) &\leq &\int_0^\infty  \int_0^\infty e^{-\alpha s}\langle f_{n_k}-  T_s^{\mu_+} f_{n_k}, f_{n_k} \rangle_m  \,d\widetilde{\nu}(\alpha) \,ds\\
\nonumber &= &\int_0^\infty \langle f_{n_k}-  T_s^{\mu_+} f_{n_k}, f_{n_k} \rangle_m  \,d\nu(s)\\
&=& \mathcal{E}^{\mu_+, \Phi}(f_{n_k},f_{n_k}) \label{eq:cri_b_2}
\end{eqnarray}

Hence, by (\ref{eq:cri_b_1}), (\ref{eq:cri_b_2}), the contraction property of the resolvent, and Remark \ref{rem_lem1}, we obtain
\begin{eqnarray}
\nonumber \lefteqn{\int_M ^\infty \int_E R_\alpha^{\mu_+}f_{n_k}R_\alpha^{\mu}f_{n_k}\,d\mu_- \,d\widetilde{\nu}(\alpha)}\\
\nonumber &\leq &\frac{\varepsilon}{1-\theta_1} \mathcal{E}^{\mu_+, \Phi}(f_{n_k},f_{n_k}) + \frac{C_\varepsilon }{1-\theta_1} \int_M ^\infty  \mathcal{E}^{\mu_+}(R_\alpha^{\mu_+}R_1^{\mu_+}f_{n_k}, R_\alpha^{\mu_+}R_1^{\mu_+}f_{n_k})\,d\widetilde{\nu}(\alpha)\\
\nonumber &\leq &\frac{\varepsilon}{1-\theta_1}  \,C_1 + \frac{C_\varepsilon }{1-\theta_1}\int_M ^\infty  \frac{ 1}{\alpha^2}\mathcal{E}^{\mu_+}(R_1^{\mu_+}f_{n_k}, R_1^{\mu_+}f_{n_k})\,d\widetilde{\nu}(\alpha)\\
\nonumber &\leq &\frac{\varepsilon}{1-\theta_1}  \,C_1 + C_\varepsilon  2C \mathcal{E}^{\mu_+, \Phi}(f_{n_k},f_{n_k}) \int_M ^\infty \frac{1}{\alpha^2}\,d\widetilde{\nu}(\alpha)\\
&\leq & \frac{\varepsilon }{1-\theta_1} \,(C_1 + 2C C_1). \label{eq:cri_b_3}
\end{eqnarray}

By (\ref{eq:critical_1_2}), the non-negativity of \(\mathcal{E}^\mu\) and Remark \ref{rem_lem1}, we have
\begin{eqnarray*}
\int_E  R_\alpha^\mu f_{n_k}\, R_\alpha^{\mu_+}f_{n_k}\,d\mu_- &\leq &\frac{1}{1-\theta_\alpha} \| R_\alpha^{\mu_+} f_{n_k} \|_{L^2(\mu_-)}^2\\
 &\leq &\frac{1}{1-\theta_\alpha} \mathcal{E}^{\mu_+}(R_\alpha^{\mu_+} f_{n_k},R_\alpha^{\mu_+} f_{n_k})\\
 &\leq &\frac{1}{1-\theta_\alpha}  \frac{C}{\alpha}\left( 1+\frac{1}{\alpha} \right) \mathcal{E}^{\mu_+, \Phi}(f_{n_k},f_{n_k}).
\end{eqnarray*}

Since \(\lim_{n\to \infty}\mathcal{E}^{\mu_+, \Phi}(f_n,f_n)=1\) and \(\theta_\alpha\) is non-increasing, it holds that, for any \(\delta \leq \alpha \leq M\),
\begin{equation}
\int_E  R_\alpha^\mu f_{n_k}\, R_\alpha^{\mu_+}f_{n_k}\,d\mu_- \leq \frac{1}{1-\theta_\delta} \frac{C}{\delta}\left( 1+\frac{1}{\delta} \right)=:C_\delta, \label{eq:cri_b_4}
\end{equation}
where \(C_\delta>0\) is depending only on \(\delta\). Since \(\widetilde{\nu}([\delta, M))\) is finite, the dominated convergence theorem yields that
\begin{equation}
\lim_{n_k \to \infty}\int_\delta ^M \int_E  R_\alpha^\mu f_{n_k}\, R_\alpha^{\mu_+}f_{n_k}\,d\mu_-\,d\widetilde{\nu}(\alpha) = 
\int_\delta ^M \int_E  R_\alpha^\mu v\, R_\alpha^{\mu_+}v\,d\mu_-\,d\widetilde{\nu}(\alpha). \label{eq:cri_b_5}
\end{equation}

Combining (\ref{eq:cri_b_3}) and (\ref{eq:cri_b_5}) with (\ref{eq:critical_1_5}), \(\mathcal{E}^{\mu_+, \Phi}(f_{n_k},f_{n_k}) - \mathcal{E}^{\mu, \Phi}(f_{n_k},f_{n_k})\) converges to \(\mathcal{E}^{\mu_+, \Phi}(v,v) - \mathcal{E}^{\mu, \Phi}(v,v)\) for the case of \(b=0\).

\medskip
Thus, in both cases, we obtain \(\mathcal{E}^{\mu_+, \Phi}(v,v) - \mathcal{E}^{\mu, \Phi}(v,v)=1\). In the rest of the proof, we treat both cases simultaneously. Since \(\mathcal{E}^{\mu_+, \Phi}(v,v) - \mathcal{E}^{\mu, \Phi}(v,v)=1\), it holds that
\[0\leq \mathcal{E}^{\mu, \Phi}(v,v) = \mathcal{E}^{\mu_+, \Phi}(v,v) - 1 \leq \varliminf_{n_k \to \infty} \mathcal{E}^{\mu_+, \Phi}(f_{n_k},f_{n_k})-1=0.\]
We remark that \(v\not = 0\) by \(\mathcal{E}^{\mu_+, \Phi}(v,v)=\mathcal{E}^{\mu_+, \Phi}(v,v)-\mathcal{E}^{\mu, \Phi}(v,v)=1.\) By the non-negativity of \(\mathcal{E}^\mu\) and the completeness of \(\Phi\), it holds that \(|v| \in D_e(\mathcal{E}^{\mu_+, \Phi}) \subset D_e(\mathcal{E}^{\mu, \Phi})\). Moreover \(\mathcal{E}^{\mu, \Phi}(|v|,|v|)=0\) holds.

We assumed \(\lambda(\mu, \Phi)=1\), so by Proposition \ref{prop:comp_spectrum}, we have 
\(\lambda(\mu)=1\). By \cite[Theorem 5.19]{T14}, \((\mathcal{E}^\mu, D(\mathcal{E}^\mu))\) is critical, that is, there exists a strictly positive function \(h\) such that \(T_t^\mu h \leq h\) holds for any \(t\geq 0\) and an \(h\)-transform \((\mathcal{E}^{\mu, h}, D(\mathcal{E}^{\mu, h}))\) is an irreducible recurrent Dirichlet form on \(L^2(E;h^2 m)\). Then \(h\) also satisfies \(T_t^{\mu, \Phi}h \leq h\) for any \(t>0\), and an \(h\)-transform \((\mathcal{E}^{\mu, \Phi, h}, D(\mathcal{E}^{\mu,\Phi, h}))\) is an irreducible Dirichlet form on \(L^2(E;h^2 m)\) by \cite[Proposition 3.4]{OS26}. By \cite[Lemma 3.8]{OS26}, \(|v|/h \in D_e(\mathcal{E}^{\mu, \Phi, h})\) and 
\[\mathcal{E}^{\mu, \Phi, h}\left(\frac{|v|}{h},\frac{|v|}{h}\right) = \mathcal{E}^{\mu, \Phi}(|v|,|v|) =0.\]
By irreducibility, there exists \(c>0\) satisfying \(|v|=ch\). Thus \(|v|\) is strictly positive, and so  \((\mathcal{E}^{\mu, \Phi}, \mathcal{D}(\mathcal{E}^{\mu, \Phi}))\) is critical.
\end{proof}

Hence we obtain the following theorem.
\begin{theorem}\label{thm:subcri_cri}
Suppose that the strong Feller property {\rm (SF)} holds for \((\mathcal{E}^{\mu_+}, D(\mathcal{E}^{\mu_+}))\). For a signed smooth Radon measure \(\mu \in \mathcal{K}_{loc}(\mathcal{E})-\mathcal{K}_{\infty}(\mathcal{E}^{\mu_+})\) and a complete Bernstein function \(\Phi\) satisfying {\rm (UC)}, the following hold.
\begin{enumerate}
\item \((\mathcal{E}^{\mu, \Phi}, \mathcal{D}(\mathcal{E}^{\mu, \Phi}))\) is subcritical if and only if \(\lambda(\mu, \Phi)>1\),
\item \((\mathcal{E}^{\mu, \Phi}, \mathcal{D}(\mathcal{E}^{\mu, \Phi}))\) is critical if and only if \(\lambda(\mu, \Phi)=1\),
\item \((\mathcal{E}^{\mu, \Phi}, \mathcal{D}(\mathcal{E}^{\mu, \Phi}))\) is supercritical if and only if \(\lambda(\mu, \Phi)<1\).
\end{enumerate}
\end{theorem}

Condition {\rm (UC)} may be difficult to verify directly. Hence, int the rest of this section, we provide sufficient conditions for (UC).
\begin{corollary}\label{cor_cri_1}
Suppose that the strong Feller property {\rm (SF)} holds for \((\mathcal{E}^{\mu_+}, D(\mathcal{E}^{\mu_+}))\). For a signed smooth Radon measure \(\mu \in \mathcal{K}_{loc}(\mathcal{E})-\mathcal{K}_{\infty}(\mathcal{E}^{\mu_+})\) and a complete Bernstein function \(\Phi\), we assume that there exists \(\beta>0\) satisfying \(\beta \|R_\beta^{\mu_+}1\|_\infty < 1\). Then {\rm (UC)} holds. In particular, if \(\lambda(\mu, \Phi)=1\), then \((\mathcal{E}^{\mu, \Phi}, \mathcal{D}(\mathcal{E}^{\mu, \Phi}))\) is critical.
\end{corollary}
Before beginning the proof, we consider an estimate for \(\theta_\alpha\).
\begin{lemma}\label{lem:theta}
There exists \(C>0\) such that, for any \(0<\alpha<1\), \((1-\theta_\alpha)^{-1} \leq C\alpha^{-1}\) holds.
\end{lemma}
\begin{proof}
For any $0<\alpha<1$ and
$f\in D(\mathcal{E}^{\mu_+})$, we obtain
\begin{eqnarray*}
\mathcal{E}_\alpha^{\mu}(f,f) &=& (1-\alpha)\mathcal{E}^{\mu}(f,f) + \alpha \mathcal{E}_1^{\mu_+}(f,f)-\alpha \int_E|f|^2\,d\mu_-\\
&\geq &\alpha \mathcal{E}_1^{\mu_+}(f,f)-\alpha \theta_1^2 \mathcal{E}_1^{\mu_+}(f,f)\\
&\geq & \alpha(1-\theta_1^2)\mathcal{E}_\alpha^{\mu_+}(f,f)
\end{eqnarray*}
and so
\[\int_E |f|^2\,d\mu_- \leq \left(1-\alpha(1-\theta_1^2)\right) \mathcal{E}_\alpha^{\mu_+}(f,f).\]
Taking the supremum over all non-zero \(f\in D(\mathcal{E}^{\mu_+})\), we obtain \(\theta_\alpha^2 \leq 1- \alpha(1-\theta_1^2).\)
Since \(0\leq\theta_\alpha<1\), it follows that
\[\frac{1}{1-\theta_\alpha} =\frac{1+\theta_\alpha}{1-\theta_\alpha^2} \leq \frac{2}{1-\theta_1^2}\, \frac{1}{\alpha}.\]
\end{proof}

\begin{proof}[Proof of Corollary \ref{cor_cri_1}]
We take \(\beta >0\) satisfying \(\beta \|R_\beta^{\mu_+}1\|_\infty < 1\). Since \(\{T_s^{\mu_+}1(x)\}_s\) is non-increasing, for any \(t>0\), we have
\[\beta R_\beta^{\mu_+}1(x) \geq \int_0^t \beta e^{-\beta s}\,T_s^{\mu_+}1(x)\,ds \geq (1-e^{-\beta t})T_t^{\mu_+}1(x).\]
We take \(t_0>0\) satisfying
\[\|T_{t_0}^{\mu_+}1\|_\infty \leq \frac{\beta\| R_\beta^{\mu_+}1\|_\infty}{1-e^{-\beta t_0}} =:c_0 <1.\]
Then we have
\begin{eqnarray*}
T_t^{\mu_+, \Phi}1(x) = \int_0^\infty T_s^{\mu_+}1(x)\,\eta_t(ds) &\leq & \int_0^{t_0} 1\,\eta_t(ds)+  \int_{t_0}^\infty c_0\,\eta_t(ds)\\
&=&1-(1-c_0)\mathbb{P}(S_t \geq t_0),
\end{eqnarray*}
where \(\{S_t\}_t\) is a conservative subordinator whose Laplace exponent is \(\Phi\) and \(\{\eta_t\}_t\) is transition semigroup of \(\{S_t\}_t\). Hence we have
\[\beta R_\beta^{\mu_+, \Phi}1(x) =\int_0^\infty  \beta e^{-\beta t}T_t^{\mu_+, \Phi}1(x)\,dt \leq 1-(1-c_0)\int_0^\infty \beta e^{-\beta t}\mathbb{P}(S_t \geq t_0)\,dt <1,\]
and so \(\beta \|R_\beta^{\mu_+, \Phi}1\|_\infty < 1\) holds.

For any \(f\in D(\mathcal{E}^{\mu_+, \Phi})\), by the Stollmann-Voigt inequality for \(\mathcal{E}_\beta^{\mu_+}\), we have
\begin{eqnarray*}
\lefteqn{\int_E |R_\alpha^{\mu_+}f|^2\,dm \leq \|R_\beta^{\mu_+}1\|_\infty\, \mathcal{E}_\beta^{\mu_+}(R_\alpha^{\mu_+}f,R_\alpha^{\mu_+}f)} \\
&=&  \|R_\beta^{\mu_+}1\|_\infty\, \mathcal{E}^{\mu_+}(R_\alpha^{\mu_+}f,R_\alpha^{\mu_+}f) + \beta \|R_\beta^{\mu_+}1\|_\infty \int_E |R_\alpha^{\mu_+}f|^2\,dm
\end{eqnarray*}
and so
\begin{equation}
\int_E |R_\alpha^{\mu_+}f|^2\,dm \leq \frac{\|R_\beta^{\mu_+}1\|_\infty }{1-\beta \|R_\beta^{\mu_+}1\|_\infty }\mathcal{E}^{\mu_+}(R_\alpha^{\mu_+}f,R_\alpha^{\mu_+}f). \label{eq:cor_1_1}
\end{equation}

Similarly, by the Stollmann-Voigt inequality for \(\mathcal{E}_\beta^{\mu_+, \Phi}\), we obtain
\begin{equation}
\int_E |f|^2\,dm \leq \frac{\|R_\beta^{\mu_+,\Phi}1\|_\infty }{1-\beta \|R_\beta^{\mu_+,\Phi}1\|_\infty }\mathcal{E}^{\mu_+,\Phi}(f,f). \label{eq:cor_1_2}
\end{equation}
By (\ref{eq:cor_1_1}) and (\ref{eq:cor_1_2}), it holds that
\begin{eqnarray*}
\lefteqn{\mathcal{E}^{\mu_+}(R_\alpha^{\mu_+}f,R_\alpha^{\mu_+}f) \leq \mathcal{E}_\alpha^{\mu_+}(R_\alpha^{\mu_+}f,R_\alpha^{\mu_+}f) = \langle f, R_\alpha^{\mu_+}f \rangle_m}\\
&\leq & \|f\|_{L^2(E;m)} \|R_\alpha^{\mu_+}f\|_{L^2(E;m)}\\
&\leq &  \sqrt{\frac{\|R_\beta^{\mu_+}1\|_\infty }{1-\beta \|R_\beta^{\mu_+}1\|_\infty }\frac{\|R_\beta^{\mu_+,\Phi}1\|_\infty }{1-\beta \|R_\beta^{\mu_+,\Phi}1\|_\infty }}\,\sqrt{\mathcal{E}^{\mu_+}(R_\alpha^{\mu_+}f,R_\alpha^{\mu_+}f)} \,\sqrt{\mathcal{E}^{\mu_+,\Phi}(f,f)}.
\end{eqnarray*}
Thus we obtain
\begin{equation}
\mathcal{E}^{\mu_+}(R_\alpha^{\mu_+}f,R_\alpha^{\mu_+}f)  \leq \frac{\|R_\beta^{\mu_+}1\|_\infty }{1-\beta \|R_\beta^{\mu_+}1\|_\infty }\frac{\|R_\beta^{\mu_+,\Phi}1\|_\infty }{1-\beta \|R_\beta^{\mu_+,\Phi}1\|_\infty }  \mathcal{E}^{\mu_+,\Phi}(f,f). \label{eq:cor_1_3}
\end{equation}

By (\ref{eq:critical_1_2}), the non-negativity for \(\mathcal{E}^\mu\), Lemma \ref{lem:theta} and (\ref{eq:cor_1_3}), for \(f \in D(\mathcal{E}^{\mu_+, \Phi})\) with \(\mathcal{E}^{\mu_+, \Phi}(f,f)\leq 1\), we have
\begin{eqnarray}
\nonumber \int_E  R_\alpha^\mu f\, R_\alpha^{\mu_+} f\,d\mu_- &\leq & \int_E \frac{1}{1-\theta_\alpha} |R_\alpha^{\mu_+} f|^2\,d\mu_-\\
\nonumber &\leq & \frac{C}{\alpha} \mathcal{E}^{\mu_+}(R_\alpha^{\mu_+}f,R_\alpha^{\mu_+}f)  \\
\nonumber &\leq & \frac{C}{\alpha} \frac{\|R_\beta^{\mu_+}1\|_\infty }{1-\beta \|R_\beta^{\mu_+}1\|_\infty }\frac{\|R_\beta^{\mu_+,\Phi}1\|_\infty }{1-\beta \|R_\beta^{\mu_+,\Phi}1\|_\infty }  \mathcal{E}^{\mu_+,\Phi}(f,f)\\
&\leq &  \frac{C}{\alpha} \frac{\|R_\beta^{\mu_+}1\|_\infty }{1-\beta \|R_\beta^{\mu_+}1\|_\infty }\frac{\|R_\beta^{\mu_+,\Phi}1\|_\infty }{1-\beta \|R_\beta^{\mu_+,\Phi}1\|_\infty }. \label{eq:cor_1_4}
\end{eqnarray}
The right-hand side of (\ref{eq:cor_1_4}) is bounded by a constant multiple of  \(\alpha^{-1}\), so, {\rm (UC)} follows from (\ref{eq:Levytilde}).
\end{proof}

\begin{corollary}\label{cor_cri_2}
Suppose that the strong Feller property {\rm (SF)} holds for \((\mathcal{E}^{\mu_+}, D(\mathcal{E}^{\mu_+}))\). For a signed smooth Radon measure \(\mu \in \mathcal{K}_{loc}(\mathcal{E})-\mathcal{K}_{\infty}(\mathcal{E}^{\mu_+})\) and a complete Bernstein function \(\Phi\), we assume that \(\|R^{\mu_+}1\|_\infty <\infty\). Then, {\rm (UC)} holds. In particular, if \(\lambda(\mu, \Phi)=1\), then \((\mathcal{E}^{\mu, \Phi}, \mathcal{D}(\mathcal{E}^{\mu, \Phi}))\) is critical.
\end{corollary}
\begin{proof}
For any \(t>0\), we have
\[R^{\mu_+}1(x) \geq \int_0^t T_s^{\mu_+}1(x)\,ds \geq t \,T_t^{\mu_+}1(x)\]
and so, we can take large \(t_0>0\) satisfying \(\|T_{t_0}^{\mu_+}1 \|_\infty \leq 1/2\). Moreover, by the Markov property, it holds that \(\|T_{nt_0}^{\mu_+}1 \|_\infty \leq 2^{-n}\) for any \(n\in \mathbb{N}\). For any \(t_0 \leq t\) satisfying \(nt_0<t\leq (n+1)t_0\), we have \(T_t^{\mu_+}1(x) \leq 2^{-n} \leq 2^{1-t/t_0}\). Then we have
\begin{eqnarray*}
R^{\mu_+, \Phi}1 (x) &=& \int_0^\infty  \int_0^\infty T_s^{\mu_+}1(x) \,\eta_t(ds) \,dt \\
&\leq & \int_0^\infty 2^{1-\frac{s}{t_0}} \,\eta_t(ds) \,dt =2\,\int_0^\infty \exp{\left(-t\Phi\left(\frac{\log 2}{t_0}\right)\right)} \,dt = 2/\Phi\left(\frac{\log 2}{t_0}\right)
\end{eqnarray*}
where \(\{\eta_t\}_t\) is the convolution semigroup of probability distributions of a subordinator whose Laplace exponent is \(\Phi\). Hence we obtain \(\|R^{\mu_+, \Phi}\|_\infty < \infty.\)

For any \(f\in D(\mathcal{E}^{\mu_+, \Phi})\), by the Stollmann-Voigt inequality for \(\mathcal{E}^{\mu_+}\), we have
\begin{equation}
\int_E |R_\alpha^{\mu_+}f|^2\,dm \leq \|R^{\mu_+}1\|_\infty\, \mathcal{E}^{\mu_+}(R_\alpha^{\mu_+}f,R_\alpha^{\mu_+}f) \label{eq:cor_2_1}
\end{equation}
and, by the Stollmann-Voigt inequality for \(\mathcal{E}^{\mu_+, \Phi}\), we obtain
\begin{equation}
\int_E |f|^2\,dm \leq \|R^{\mu_+,\Phi}1\|_\infty \mathcal{E}^{\mu_+,\Phi}(f,f). \label{eq:cor_2_2}
\end{equation}
By (\ref{eq:cor_2_1}) and (\ref{eq:cor_2_2}), it holds that
\begin{eqnarray*}
\lefteqn{\mathcal{E}^{\mu_+}(R_\alpha^{\mu_+}f,R_\alpha^{\mu_+}f) \leq \mathcal{E}_\alpha^{\mu_+}(R_\alpha^{\mu_+}f,R_\alpha^{\mu_+}f) = \langle f, R_\alpha^{\mu_+}f \rangle_m}\\
&\leq & \|f\|_{L^2(E;m)} \|R_\alpha^{\mu_+}f\|_{L^2(E;m)}\\
&\leq &  \sqrt{\|R^{\mu_+}1\|_\infty  \|R^{\mu_+,\Phi}1\|_\infty}\,\sqrt{\mathcal{E}^{\mu_+}(R_\alpha^{\mu_+}f,R_\alpha^{\mu_+}f)} \,\sqrt{\mathcal{E}^{\mu_+,\Phi}(f,f)}
\end{eqnarray*}
and it holds that
\begin{equation}
\mathcal{E}^{\mu_+}(R_\alpha^{\mu_+}f,R_\alpha^{\mu_+}f)  \leq  \|R^{\mu_+}1\|_\infty  \|R^{\mu_+,\Phi}1\|_\infty\, \mathcal{E}^{\mu_+,\Phi}(f,f). \label{eq:cor_2_3}
\end{equation}
Combining this with Lemma \ref{lem:theta}, we have {\rm (UC)}.
\end{proof}

A Dirichlet form \((\mathcal{E}^{\mu_+}, D(\mathcal{E}^{\mu_+}))\) is said to be in class (T), if it is irreducible, strong Feller and satisfies the following tightness property:

For any \(\varepsilon >0\), there exists a compact set \(K\) such that \(\sup_{x\in E}R_1^{\mu_+}1_{K^c}(x) \leq \varepsilon\).
\begin{corollary}\label{cor_cri_3}
Suppose that \((\mathcal{E}^{\mu_+}, D(\mathcal{E}^{\mu_+}))\) is in class {\rm (T)} for a signed smooth Radon measure \(\mu \in \mathcal{K}_{loc}(\mathcal{E})-\mathcal{K}_{\infty}(\mathcal{E}^{\mu_+})\). Then {\rm (UC)} holds. In particular, if \(\lambda(\mu, \Phi)=1\), then \((\mathcal{E}^{\mu, \Phi}, \mathcal{D}(\mathcal{E}^{\mu, \Phi}))\) is critical.
\end{corollary}
\begin{proof}
We prove that there exists \(C>0\) such that, for any \(f\in D(\mathcal{E}^{\mu_+})\), it holds that 
\begin{equation}
\|f\|^2_{L^2(E;m)} \leq C \mathcal{E}^{\mu_+}(f,f). \label{eq:cor_3_1}
\end{equation}
If this does not hold, there exists \(f_n \in D(\mathcal{E}^{\mu_+})\) satisfying \(\|f_n\|_{L^2(E;m)}=1\) and \(\mathcal{E}^{\mu_+}(f_n,f_n) \to 0\). Then \(\{\mathcal{E}_1^{\mu_+}(f_n,f_n)\}_n\) is bounded and so we can take a subsequence \(\{n_k\}\) such that \(f_{n_k}\) converges weakly to some \(u \in D(\mathcal{E}^{\mu_+})\) in \(\mathcal{E}_1^{\mu_+}\).  Since \((\mathcal{E}^{\mu_+}, D(\mathcal{E}^{\mu_+}))\) belongs to class {\rm (T)}, by \cite[Corollary 4.5]{T19}, \((D(\mathcal{E}^{\mu_+}),\mathcal{E}_1^{\mu_+})\) is compactly embedded into \(L^2(E;m)\). Then \(f_{n_k}\) converges to \(u\) in \(L^2(E;m)\), and so \(\|u\|_{L^2(E;m)} = 1\). However, it holds that 
\[0\leq \mathcal{E}^{\mu_+}(u,u) \leq \varliminf_{n_k \to \infty} \mathcal{E}^{\mu_+}(f_{n_k},f_{n_k}) = 0\]
and, by transience of \((\mathcal{E}^{\mu_+}, D(\mathcal{E}^{\mu_+}))\) , it holds that \(u=0\). This is a contradiction, so (\ref{eq:cor_3_1}) holds.

If \(b>0\), we also obtain 
\begin{equation}
\|f\|^2_{L^2(E;m)} \leq C \mathcal{E}^{\mu_+}(f,f) \leq \frac{C}{b} \mathcal{E}^{\mu_+, \Phi}(f,f). \label{eq:cor_3_2}
\end{equation}

By the spectral decomposition and (\ref{eq:cor_3_1}), we have
\begin{eqnarray*}
\frac{d}{dt} \|T_t^{\mu_+} f\|_{L^2(E;m)}^2 = -2\,\mathcal{E}^{\mu_+}(T_t^{\mu_+}f,T_t^{\mu_+}f) \leq -\frac{2}{C}\,\|T_t^{\mu_+}f\|^2_{L^2(E;m)}
\end{eqnarray*}
and, by Gronwall's inequality, we obtain \(\|T_t^{\mu_+}f \|_{L^2(E;m)} \leq e^{-t/C}\,\|f\|_{L^2(E;m)}\) for any \(t>0\). By the Cauchy--Schwarz inequality, we have
\[\langle f-T_t^{\mu_+}f, f\rangle_m \geq (1-e^{-t/C})\,\|f\|_{L^2(E;m)}^2.\]
Hence, if \(b=0\), by \cite[Theorem 2.1]{O02}, we have
\begin{equation}
\mathcal{E}^{\mu_+, \Phi}(f,f) = \int_0^\infty \langle f-T_t^{\mu_+}f, f\rangle_m \,d\nu(t) \geq \int_0^\infty \langle(1-e^{-t/C})\,d\nu(t) \|f\|_{L^2(E;m)}^2 = \Phi(C^{-1}) \|f\|_{L^2(E;m)}^2. \label{eq:cor_3_3}
\end{equation}
By (\ref{eq:cor_3_2}) and (\ref{eq:cor_3_3}), for any complete Bernstein function \(\Phi\), there exists \(C_\Phi > 0\) such that, for any \(f \in D(\mathcal{E}^{\mu_+})\), it holds that 
\begin{equation}
\|f\|^2_{L^2(E;m)} \leq C_\Phi \mathcal{E}^{\mu_+, \Phi}(f,f). \label{eq:cor_3_4}
\end{equation}

Similarly to the proof of (\ref{eq:cor_2_3}), by (\ref{eq:cor_3_1}) and (\ref{eq:cor_3_4}), we have
\[
\mathcal{E}^{\mu_+}(R_\alpha^{\mu_+}f,R_\alpha^{\mu_+}f)  \leq  CC_\Phi\, \mathcal{E}^{\mu_+,\Phi}(f,f).
\]
Combining this with Lemma \ref{lem:theta}, we have {\rm (UC)}.
\end{proof}

Since, for \(\mu_+\in \mathcal{K}(\mathcal{E})\), the strong Feller property for \((\mathcal{E}^{\mu_+}, D(\mathcal{E}^{\mu_+}))\) follows from the strong Feller property for \((\mathcal{E}, D(\mathcal{E}))\) (\cite[Corollary 1.2]{CK09}), the following corollary holds.
\begin{corollary}\label{cor_cri_4}
Suppose that \((\mathcal{E}, D(\mathcal{E}))\) is in class {\rm (T)}. Then, for a signed smooth Radon measure \(\mu \in \mathcal{K}(\mathcal{E})-\mathcal{K}_{\infty}(\mathcal{E}^{\mu_+})\), {\rm (UC)} holds. In particular, if \(\lambda(\mu, \Phi)=1\), then \((\mathcal{E}^{\mu, \Phi}, \mathcal{D}(\mathcal{E}^{\mu, \Phi}))\) is critical.
\end{corollary}

\section{Example}\label{sec:Example}
In this section, we compute \(\lambda(\mu, \Phi)\) for a subordinated Schr\"{o}dinger operators. The first example shows that the necessary and sufficient condition for a complete Bernstein function that preserving \(\lambda(\mu, \Phi)=1\) when the original Schr\"{o}dinger operators corresponds to the critical (fractional) Hardy inequality on \(\mathbb{R}^d\). The second example shows that \(\lambda(\mu, \Phi)=1\) for a subordinated Laplacian with a bounded integrable potential on \(\mathbb{R}^d\) for \(d\geq 5\) by stable subordinators.

\begin{example}[The Hardy inequality for \((-\Delta)^{\alpha/2}\)]\label{ex.Hardy}
We consider the case corresponding to the critical Hardy inequality for \(\mathcal{L}:=-\frac{1}{2}(-\Delta)^{\alpha/2}\) on \(\mathbb{R}^d\) with \(0<\alpha \leq 2\) and \(\alpha< d\). A Dirichlet form \((\mathcal{E}, \mathcal{D}(\mathcal{E}))\) on \(L^2(\mathbb{R}^d ;dx)\) associated with \(\mathcal{L}:=-\frac{1}{2}(-\Delta)^{\alpha/2}\) is denoted by, for \(\alpha=2\),
\begin{eqnarray*}
\begin{cases}
\mathcal{E}(f,g):= \frac{1}{2} \int_{\mathbb{R}^d} \nabla f\cdot \nabla g \,dx \\
\mathcal{D}(\mathcal{E}):=H^1(\mathbb{R}^d):=\{f\in L^2(\mathbb{R}^d) : \nabla f \in L^2(\mathbb{R}^d)\}
\end{cases}
\end{eqnarray*}
where the derivatives are taken in the Schwartz distribution sense, and, for \(0<\alpha <2\),
\begin{eqnarray*}
\begin{cases}
\mathcal{E}(f,f):= \frac{\mathcal{A}(d,\alpha)}{2}\iint_{\mathbb{R}^d\times \mathbb{R}^d \setminus \rm{diag}} \frac{(f(x)-f(y))^2}{|x-y|^{d+\alpha}}\,dx\,dy,\\
\mathcal{D}(\mathcal{E}) := \left\{f\in L^2(\mathbb{R}^d) : \iint_{\mathbb{R}^d\times \mathbb{R}^d \setminus \rm{diag}} \frac{(f(x)-f(y))^2}{|x-y|^{d+\alpha}}\,dx\,dy<\infty \right\},
\end{cases}
\end{eqnarray*}
where \({\rm diag}:=\{(x,y)\in \mathbb{R}^d\times \mathbb{R}^d : x=y\}\) and
\[\mathcal{A}(d,\alpha):= \frac{\alpha 2^{\alpha-1} \Gamma((d+\alpha)/2)}{\pi^{d/2} \Gamma(1-\alpha/2)}.\]

In both cases, \((\mathcal{E}, \mathcal{D}(\mathcal{E}))\) is an irreducible transient regular Dirichlet form satisfying the strong Feller property (SF). It is  well-known that the following Hardy inequality holds. 
\begin{eqnarray}
\lambda_* \int_{\mathbb{R}^d} \frac{|f(x)|^2}{|x|^\alpha}\,dx \leq \mathcal{E}(f,f) \label{eq:ex_Hardy_1}
\end{eqnarray}
for \(f\in \mathcal{D}(\mathcal{E})\), where 
\[\lambda_*:=2^{\alpha-1}\left(\frac{\Gamma(\frac{d+\alpha}{4})}{\Gamma(\frac{d-\alpha}{4})}\right)^2.\]

Let \(d\mu_\lambda(x) := \lambda |x|^{-\alpha}\,dx\). Then a Schr\"{o}dinger form \((\mathcal{E}^{-\mu_{\lambda_*}}, \mathcal{D}(\mathcal{E}^{-\mu_{\lambda_*}}))\) is critical. For \(\lambda <\lambda_*\) (resp. \(\lambda >\lambda_*\)), a Schr\"{o}dinger form \((\mathcal{E}^{-\mu_\lambda}, \mathcal{D}(\mathcal{E}^{-\mu_\lambda}))\) is subcritical (resp. supercritical).

We consider \(\lambda(\mu, \Phi)\) for the case of \(\mu_+:=0\), \(\mu:=-\mu_-:=-\mu_{\lambda_*}\) and a complete Bernstein function \(\Phi\). We note that \(\lambda(\mu)=1\) in this case. Then the following holds.

\begin{proposition}\label{prop:ex_1}
For a complete Bernstein function \(\Phi\), \(\lambda(\mu, \Phi)>1\) if and only if \(\sup_{s>0}\limits \frac{s\Phi'(s)}{\Phi(s)}<1\).
\end{proposition}
\begin{proof}
First we show that, a complete Bernstein function \(\Phi\) satisfies  \(\sup_{s>0}\limits \frac{s\Phi'(s)}{\Phi(s)}<1\) if and only if, there exist \(C\geq 1\) and \(0<\beta <2\) such that \(\Phi(Rs)\leq CR^{\beta/2}\Phi(s)\) holds for any \(R\geq 1\).  If \(\sup_{s>0}\limits \frac{s\Phi'(s)}{\Phi(s)}\leq \beta/2\) holds for \(0<\beta <2\), then we have \((\log{\Phi(s)})' \leq \beta/2s\) and, by integrating from \(s\) to \(Rs\), it holds that \(\Phi(Rs)\leq R^{\beta/2}\Phi(s)\). On the other hand, suppose that there exist \(C\geq 1\) and \(0<\beta <2\) such that \(\Phi(Rs)\leq CR^{\beta/2}\Phi(s)\) holds for any \(R\geq 1\). Since \(\Phi\) is concave, we have \(\Phi(s)-\Phi(s/L)\geq (s-s/L)\Phi'(s)\) for \(L>1\). Then we have
\[\frac{s\Phi'(s)}{\Phi(s)} \leq \frac{1-L^{-\beta/2}C^{-1}}{1-L^{-1}},\]
and by taking \(L>1\) satisfying \(L^{1-\beta/2}>C,\)  \(\sup_{s>0}\limits \frac{s\Phi'(s)}{\Phi(s)}<1\) holds.

We assume that \(\sup_{s>0}\limits \frac{s\Phi'(s)}{\Phi(s)}<1\) holds. Then we can take   \(C\geq 1\) and \(0<\beta <2\) such that \(\Phi(Rs)\leq CR^{\beta/2}\Phi(s)\) holds for any \(R\geq 1\). We define \(\Psi\) by
\[\Psi(\lambda):=\int_0^\infty \frac{\lambda}{\lambda + r^{\theta}}\Phi'(r)\,dr=\int_0^\infty \frac{\lambda}{\lambda + t}\,d\Phi(t^\frac{1}{\theta})\]
where \(\theta := \beta/4 + 1/2\). Then we have \(\Phi(t^{1/\theta}) \leq C t^{\beta/(2\theta)}\Phi(1)\) and
\begin{eqnarray*}
\int_0^\infty \frac{d\Phi(t^\frac{1}{\theta})}{1+t} &\leq & \Phi(1) + \int_1^\infty \frac{d\Phi(t^\frac{1}{\theta})}{t}\\
&=& \Phi(1) +\left[\frac{\Phi(t^\frac{1}{\theta})}{t} \right]_1^\infty +  \int_1^\infty \frac{\Phi(t^\frac{1}{\theta})}{t^2}\,dt\\
&\leq & C\Phi(1) \, \int_1^\infty t^{\frac{\beta}{2\theta}-2}\,dt <\infty.
\end{eqnarray*}
By \cite[Theorem 6.2]{SSV12}, \(\Psi\) is also a complete Bernstein function. Moreover we have
\[\Psi(\lambda) \leq  \lambda\,\int_\lambda^\infty \frac{\Phi(t^\frac{1}{\theta})}{t^2}\,dt \leq C \lambda \Phi(\lambda^{\frac{1}{\theta}})\int_\lambda^\infty \left(\frac{r}{\lambda} \right)^{\frac{\beta}{2\theta}}\frac{dr}{r^2} = \frac{C}{1-\frac{\beta}{2\theta}}\Phi(\lambda^{\frac{1}{\theta}})\]
and
\[\Psi(\lambda)\geq \int_0^\lambda  \frac{\lambda}{\lambda + t}\,d\Phi(t^\frac{1}{\theta}) \geq \frac{1}{2}\Phi(\lambda^{\frac{1}{\theta}}).\]

By \cite[Theorem 1.2]{KMVZZ18} and \cite[Theorem 1.1]{FMS21}, there exists \(c>1\) satisfying \(\mathcal{E}^{ s^{\theta}}(f,f) \leq c \,\mathcal{E}^{\mu, s^{\theta}}(f,f)\) holds for any \(f\in C_c^\infty.\) We remark that \(\mathcal{E}^{ s^{\theta}}(f,f)\) corresponds to \((-\Delta)^{\alpha \theta/2}\) and \(\mathcal{E}^{\mu, s^{\theta}}(f,f)\) corresponds to \((-(-\Delta)^{\alpha/2}-\mu)^{\theta}\). Since \(\Psi\) is complete and concave, \(\mathcal{E}^{ \Psi(s^{\theta})}(f,f) \leq c \,\mathcal{E}^{\mu, \Psi(s^{\theta})}(f,f)\) also holds. Since \(\Psi(\lambda^\theta)\) is comparable to \(\Phi(\lambda)\), there exists \(c_1>1\) satisfying \(\mathcal{E}^{\Phi}(f,f) \leq c_1 \,\mathcal{E}^{\mu, \Phi}(f,f)\) holds for any \(f\in C_c.\) Hence we have
\[\lambda(\mu, \Phi)=\inf_{f\in D(\mathcal{E})\cap C_c}\frac{\mathcal{E}^\Phi(f,f)}{\mathcal{E}^\Phi(f,f)-\mathcal{E}^{\mu,\Phi}(f,f)} \geq \frac{1}{1-c_1^{-1}}>1.\]

Next we assume that \(\sup_{s>0}\limits \frac{s\Phi'(s)}{\Phi(s)}\geq 1\) holds. We note that, since \(\Phi\) is non-decreasing and concave,  \(\sup_{s>0}\limits \frac{s\Phi'(s)}{\Phi(s)}= 1\) holds. Then there exists a sequence \(\{s_n\}_n\) satisfying \(1-n^{-2} \leq s_n\Phi'(s_n)/\Phi(s_n) \leq 1\). We fix \(n \in \mathbb{N}\). For the convenience, we set \(R:=(s_n/n)^{1/\alpha}\). We define 
\[\Phi_n(t):= \frac{\Phi\left(\frac{s_n}{n}t\right)}{\Phi\left(\frac{s_n}{n}\right)} = \frac{\Phi\left(R^\alpha t\right)}{\Phi\left(R^\alpha\right)},\]
then \(\Phi_n\) is also a complete Bernstein function. Since \(\Phi\) is concave, for any \(t>0\) and \(n >t\), we have
\[\Phi\left(s_n \frac{t}{n}\right) \geq \frac{t}{n} \Phi(s_n) + \left(1-\frac{t}{n}\right)\Phi(0) =  \frac{t}{n} \Phi(s_n)\]
and by \(1-n^{-2} \leq s_n\Phi'(s_n)\), 
\[\Phi\left(s_n \frac{t}{n}\right) \leq  \Phi(s_n) + \left(\frac{t}{n}-1\right)s_n\Phi'(s_n) \leq \left(\frac{t}{n}+\frac{1}{n^2}\right) \Phi(s_n).\]
By taking \(t=1\), we have \(\Phi(s_n)/n \leq \Phi(s_n/n) \leq (n^{-1}+n^{-2})\Phi(s_n).\) Hence we have \(t/(1+n^{-1}) \leq \Phi_n(t) \leq t + 1/n\) and so \(\Phi_n\) converges to {\rm id}.

For \(f\in D(\mathcal{E})\cap C_c\), we set 
\[f_n(x):= \left(\frac{s_n}{n}\right)^{\frac{d}{2\alpha}} f\left(\left(\frac{s_n}{n}\right)^{\frac{1}{\alpha}}x\right) = R^{d/2}f(Rx).\]
Then we have
\[\mathcal{E}^{\mu}(f_n,f_n) = \frac{s_n}{n} \,\mathcal{E}^\mu(f,f)= R^\alpha\,\mathcal{E}^\mu(f,f)\]
and the correspondence between \((\mathcal{E}^\mu, D(\mathcal{E}^\mu))\) and \(\{T^\mu_t\}_t\) yields that  
\[T_s^\mu f_n(x) =  \left(\frac{s_n}{n}\right)^{\frac{d}{2\alpha}}\,T_{ \frac{s_n}{n}   s}^\mu f\left(\left(\frac{s_n}{n}\right)^{\frac{1}{\alpha}}x\right) = R^{d/2}\,T_{R^\alpha s}^\mu f(Rx).\]
By \cite[Theorem 3.2]{OS26}, we have
\begin{eqnarray*}
\mathcal{E}^{\mu, \Phi}(f_n,f_n) &=& b\,\mathcal{E}^{\mu}(f_n,f_n) + \int_0^\infty \langle f_n-T_s^\mu f_n, f_n\rangle_m\,d\nu(s)\\
&=&bR^\alpha\,\mathcal{E}^{\mu}(f,f) + R^d\,\int_0^\infty \langle f(R\cdot)-T_{R^\alpha s}^\mu f(R\cdot), f(R\cdot)\rangle_m\,d\nu(s)\\
&=&bR^\alpha\,\mathcal{E}^{\mu}(f,f) +  \,\int_0^\infty \langle f-T_{R^\alpha s}^\mu f, f\rangle_m\,d\nu(s)\\
&=& \Phi(R^\alpha)\,\mathcal{E}^{\mu, \Phi_n}(f,f).
\end{eqnarray*}
Similarly, \(\mathcal{E}^{\Phi}(f_n,f_n)=\Phi(R^\alpha)\,\mathcal{E}^{\Phi_n}(f,f)\) holds. Then we have
\begin{equation}
\frac{\mathcal{E}^{\Phi}(f_n,f_n)}{\mathcal{E}^{\Phi}(f_n,f_n)-\mathcal{E}^{\mu, \Phi}(f_n,f_n)} = \frac{\mathcal{E}^{\Phi_n}(f,f)}{\mathcal{E}^{\Phi_n}(f,f)-\mathcal{E}^{\mu, \Phi_n}(f,f)} \label{eq:ex_1_1}
\end{equation}
Since \(\lambda(\mu)=1\), we can take \(h^{(k)}\in D(\mathcal{E})\cap C_c\) satisfying 
\[\lim_{k\to \infty} \frac{\mathcal{E}(h^{(k)},h^{(k)})}{\mathcal{E}(h^{(k)},h^{(k)})-\mathcal{E}^\mu(h^{(k)},h^{(k)})} = 1.\]
Since \(\Phi_n\) converges to {\rm id}, by (\ref{eq:ex_1_1}), we can take a subsequence \(h^{(n_k)}\) satisfying
\[\lim_{n_k \to \infty} \frac{\mathcal{E}^{\Phi}(h^{(n_k)},h^{(n_k)})}{\mathcal{E}^{\Phi}(h^{(n_k)},h^{(n_k)})-\mathcal{E}^{\mu, \Phi}(h^{(n_k)},h^{(n_k)})} = 1\]
and then \(\lambda(\mu, \Phi) = 1.\)
\end{proof}

For example, a stable subordinator \(\Phi(s)=s^{\beta/2}\) with \(0<\beta <2\) satisfies \(\frac{s\Phi'(s)}{\Phi(s)}=\beta/2 <1\), so \((-\frac{1}{2}(-\Delta)^{\alpha/2}+\mu)^{\beta/2}\) is subcritical. For other examples satisfying the condition \(\sup_{s>0}\limits \frac{s\Phi'(s)}{\Phi(s)}<1\), there are \(\Phi(s) = \log(1+s^{\beta/2})\), \(\Phi(s) = s^{\beta/2}/(1+s^{\beta/2})\), \(\Phi(s) = s^{\beta/2}+s^{\gamma/2}\), and \(\Phi(s) = (m+s^{\beta/2})^{\gamma/2}-m^{\gamma/2}\) for \(0<\beta, \gamma <2\) and \(m>0\). Hence, for these \(\Phi\),  \(\Phi(-\frac{1}{2}(-\Delta)^{\alpha/2}+\mu)\) are subcritical.

For a complete Bernstein function \(\Phi(s)= \log(1+s)\), \(\Phi(s)= s/(1+s)\),  \(\Phi(s) = s+s^{\gamma/2}\), and \(\Phi(s) = (m+s)^{\gamma/2}-m^{\gamma/2}\)  for \(0< \gamma <2\) and \(m>0\), it holds that \(\sup_s \frac{s\Phi'(s)}{\Phi(s)} =1\) and so 
\(\lambda(\mu, \Phi) = 1\) holds. In particular, for a complete Bernstein function \(\Phi\) with a positive drift \(b>0\) and a L\'{e}vy measure \(\nu\), we have 
\[\frac{s\Phi'(s)}{\Phi(s)} \geq \frac{b}{b + \int_0^\infty \frac{1-e^{-st}}{s}\,d\nu(t)} \to 1\]
as \(s \to \infty\), so \(\lambda(\mu, \Phi) = 1\) holds.
\end{example}

\begin{example}\label{ex:2}
We consider the case corresponding to the critical Hardy inequality for \(\mathcal{L}:=\frac{1}{2}\Delta\) on \(\mathbb{R}^d\) with \(3\leq d\). A Dirichlet form \((\mathcal{E}, \mathcal{D}(\mathcal{E}))\) on \(L^2(\mathbb{R}^d ;dx)\) associated with \(\mathcal{L}:=\Delta\) is denoted by
\begin{eqnarray*}
\begin{cases}
\mathcal{E}(f,g):= \frac{1}{2} \int_{\mathbb{R}^d} \nabla f\cdot \nabla g \,dx \\
\mathcal{D}(\mathcal{E}):=H^1(\mathbb{R}^d):=\{f\in L^2(\mathbb{R}^d) : \nabla f \in L^2(\mathbb{R}^d)\}
\end{cases}
\end{eqnarray*}
where the derivatives are taken in the Schwartz distribution sense. A Dirichlet form \((\mathcal{E}, \mathcal{D}(\mathcal{E}))\) is irreducible, transient and regular. We may assume \(V\neq 0\). For a non-negative bounded function \(V\in L^1(\mathbb{R}^d)\), we set \(d\nu:=-d\nu_-:=-V\,dx.\) By the Sobolev inequality,
\[0<\lambda(\nu)=\inf_{f\in D(\mathcal{E})\cap C_c} \frac{\mathcal{E}(f,f)}{\int |f|^2 V\,dx} < \infty.\]
Hence, by setting \(\mu:=\lambda(\nu)\nu\), we have \(\lambda(\mu)=1\). Since \(\mu \in \mathcal{K}_\infty\) (See \cite[Example 6.1]{KM25} for example), we have \(D(\mathcal{E}) = D(\mathcal{E}^\mu)\) and, by \cite[Theorem 5.19]{T14}, \(\frac{1}{2}\Delta - \mu= \frac{1}{2}\Delta - \lambda(Vm)V\) is critical. We take a strictly positive function \(h \in D_e(\mathcal{E}^\mu) = D_e(\mathcal{E})\) satisfying \(\mathcal{E}^\mu(h,h)=0\). In particular, \(Vh^2 \in L^1(\mathbb{R}^d)\) holds. By the Cauchy--Schwarz inequality, \(Vh\in L^1(\mathbb{R}^d)\cap L^2(\mathbb{R}^d)\) holds, and so its Fourier transform \(\widehat{Vh} \in L^\infty(\mathbb{R}^d)\cap L^2(\mathbb{R}^d)\). By the Sobolev inequality, it holds that \(h\in L^{2d/(d-2)}\). Hence we can consider the Fourier transform \(\widehat{h}\) of \(h\). By \(\mathcal{E}^\mu(h,h)=0\) and the Fourier transform, we obtain
\[\widehat h(\xi) = \frac{2\lambda(\nu)}{|\xi|^2} \widehat{Vh}(\xi).\]

For \(n \in \mathbb{N}\), we choose \(f_n\in D(\mathcal{E})\) satisfying \(\widehat{f_n}(\xi)=\mathbf{1}_{\{|\xi|\geq1/n\}}\widehat{h}(\xi)\) where \(\widehat{f_n}\) is the Fourier transform of \(f_n\). Indeed, by Plancherel's theorem,
\[\int |f_n|^2\,dm = 4\lambda(\nu)^2 \int_{|\xi|\geq1/n} \frac{|\widehat{Vh}(\xi)|^2}{|\xi|^4}\,d\xi \leq C\int_{1/n}^1r^{d-5}\,dr +C\int_{|\xi|\geq1}|\widehat{Vh}(\xi)|^2\,d\xi < \infty\]
and
\[\int_{\mathbb R^d}|\nabla f_n|^2\,dx = \int_{\mathbb R^d}|\xi|^2|\widehat{f_n}(\xi)|^2\,d\xi = \int_{|\xi|\geq1/n}|\xi|^2|\widehat h(\xi)|^2\,d\xi \leq  \int_{\mathbb R^d}|\nabla h|^2\,dx < \infty.\]

By the non-negativity of \(\mathcal{E}^\mu\) and the identity \(\mathcal{E}^\mu(h,f)=0\), we have
\begin{eqnarray*}
\mathcal{E}^\mu(f_n,f_n) &=& \mathcal{E}^\mu(h-f_n,h-f_n) \leq  \mathcal{E}(h-f_n,h-f_n)\ =\  \frac12\int_{|\xi|<1/n} |\xi|^2|\widehat h(\xi)|^2\,d\xi\\ 
&=& 2\lambda(\nu)^2 \int_{|\xi|<1/n}\frac{|\widehat{Vh}(\xi)|^2}{|\xi|^2}\,d\xi\ \leq C \int_{|\xi|<1/n}\frac{1}{|\xi|^2}\,d\xi\ \leq \  Cn^{-(d-2)}.
\end{eqnarray*}
Moreover,

\begin{eqnarray*}
\|f_n\|_{L^2}^2 = 4\lambda(\nu)^2 \int_{|\xi|\geq 1/n} \frac{|\widehat{Vh}(\xi)|^2}{|\xi|^4}\,d\xi \leq C\begin{cases}
n,&d=3,\\
1+\log n,&d=4,\\
1,&d\geq5.
\end{cases}
\end{eqnarray*}

Since \(s^{\beta/2}\) is concave for \(0<\beta<2\), \(s^{\beta/2} \leq \frac{\beta}{2}t^{\beta/2-1}s + (1-\frac{\beta}{2} )t^{\beta/2}\) holds for any \(s\geq 0\) and \(t>0\). by the spectral decomposition, we obtain
\[
\mathcal{E}^{\mu,s^{\beta/2}}(f_n,f_n) \leq \frac{\beta}{2}t^{\beta/2-1}
\mathcal{E}^\mu(f_n,f_n)+ \left(1-\frac{\beta}{2}\right)t^{\beta/2}
\|f_n\|_{L^2}^2.\]
By taking \(t=\frac{\mathcal{E}^\mu(f_n,f_n)}{\|f_n\|_{L^2}^2}\), we obtain
\[\mathcal{E}^{\mu,s^{\beta/2}}(f_n,f_n) \leq \mathcal{E}^\mu(f_n,f_n)^{\beta/2} \|f_n\|_{L^2}^{2-\beta}. \]
Consequently, it holds that
\[\mathcal{E}^{\mu,\Phi}(f_n,f_n) \leq C\begin{cases}
n^{1-\beta},&d=3,\\
n^{-\beta}(1+\log n)^{1-\beta/2},&d=4,\\
n^{-(d-2)\beta/2},&d\geq5.
\end{cases}
\]

On the other hand, by the Fourier representation of the subordinated Dirichlet form for \(\beta\)-stable processes, it holds that
\begin{eqnarray*}
\mathcal{E}^{s^{\beta/2}}(f_n,f_n) &=& 2^{-\beta/2} \int_{\mathbb{R}^d}|\xi|^\beta
|\widehat{f_n}(\xi)|^2\,d\xi\\
&=&2^{2-\beta/2}\lambda(\nu)^2 \int_{1/n \leq |\xi| \leq 1}|\xi|^{\beta-4}|\widehat{Vh}(\xi)|^2\,d\xi.
\end{eqnarray*}
We assume that \(1<\beta<2\) for \(d=3\) and, \(0<\beta <2\) for \(d\geq4\). Since \(h\) is strictly positive, there exists a constant \(C>0\) such that \(\mathcal{E}^{s^{\beta/2}}(f_n,f_n)\geq C.\)
Therefore, we obtain
\begin{equation}
\frac{\mathcal{E}^{ s^{\beta/2}}(f_n,f_n)}
{\mathcal{E}^{s^{\beta/2}}(f_n,f_n)-\mathcal{E}^{\mu,s^{\beta/2}}(f_n,f_n)}
\longrightarrow 1\qquad (n\to\infty). \label{eq:ex_2_a}
\end{equation}
For the case \(\beta = 1\) and \(d=3\), \(\mathcal{E}^{s^{1/2}}(f_n,f_n) \geq \log n -C\) holds, and so (\ref{eq:ex_2_a}) holds for \(\beta = 1\) and \(d=3\). Hence we obtain \(\lambda(\mu, s^{\beta/2}) = 1\) for \(1\leq \beta <2\) when \(d=3\) and \(0<\beta <2\) when \(d\geq 4\). 

Next we check {\rm (UC)} for \(d\geq 5\). For \(\Phi(s)=s^{\beta/2}\), we have
\[d\widetilde{\nu}(\alpha)=\frac{\sin{(\pi \beta/2)}}{\pi}\alpha^{\beta/2}\ d\alpha.\]
We take \(0<\alpha < \delta<1\). For \(f\in D(\mathcal{E}^{s^{\beta/2}})\) with \(\mathcal{E}^{s^{\beta/2}}(f,f)\leq 1\), we set \(R_\alpha^{(1)} f(x)=\int_0^1 e^{-\alpha t}T_tf(x)\,dt\) and \(R_\alpha^{(2)} f(x)=\int_1^\infty e^{-\alpha t}T_tf(x)\,dt\). Since \(\lambda e^{-2t\lambda} \leq C_\beta t^{\beta/2-1} \lambda^{\beta/2}\), \(\sqrt{\mathcal{E}(T_tf, T_tf)} \leq C_\beta t^{\beta/4-1/2} \sqrt{\mathcal{E}^{s^{\beta/2}}(f,f)}\) follows from the spectral decomposition. By the Sobolev inequality and the spectral decomposition again, we have
\begin{eqnarray*}
\sqrt{\int_{\mathbb R^d}|R_\alpha^{(1)} f|^2V\,dx} &\leq &\|V\|_{L^{d/2}}^{1/2} \|R_\alpha^{(1)} f\|_{L^{2d/(d-2)}} \\
&\leq &C\|V\|_{L^{d/2}}^{1/2}\sqrt{\mathcal{E}(R_\alpha^{(1)} f,R_\alpha^{(1)} f)} \\
&\leq &C\|V\|_{L^{d/2}}^{1/2} \int_0^1e^{-\alpha t}\sqrt{\mathcal{E}(T_tf,T_tf)}\,dt \\
&\leq &C_\beta \int_0^1 e^{-\alpha t}t^{\beta/4-1/2}\,dt \,\sqrt{\mathcal{E}^{s^{\beta/2}}(f,f)}\\
&\leq & C_\beta.
\end{eqnarray*}
By using a heat kernel estimate of Brownian motion and the Sobolev inequality, we have 
\[\|T_tf\|_\infty \leq Ct^{-(d-\beta)/4} \|f\|_{L^{2d/(d-\beta)}}\leq  Ct^{-(d-\beta)/4} \sqrt{\mathcal{E}^{s^{\beta/2}}(f,f)}.\]
Hence we have
\begin{eqnarray*}
\sqrt{\int_{\mathbb R^d}|R_\alpha^{(2)} f|^2V\,dx} &\leq &\|V\|_{L^{1}}^{1/2} \|R_\alpha^{(2)} f\|_{\infty} \\
&\leq & C\|V\|_{L^{1}}^{1/2} \int_1^\infty e^{-\alpha t} t^{-(d-\beta)/4}\,dt\, \sqrt{\mathcal{E}^{s^{\beta/2}}(f,f)}\\
&\leq & C \int_1^\infty e^{-\alpha t} t^{-(d-\beta)/4}\,dt.
\end{eqnarray*}
By calculation, we obtain
\begin{eqnarray*}
\int_1^\infty e^{-\alpha t} t^{-(d-\beta)/4}\,dt &\leq & \int_1^{\alpha^{-1}}  t^{-(d-\beta)/4}\,dt + \int_{\alpha^{-1}}^\infty e^{-\alpha t} t^{-(d-\beta)/4}\,dt \\
&\leq &\int_1^{\alpha^{-1}}  t^{-(d-\beta)/4}\,dt + \alpha^{(d-\beta)/4-1} \int_1^\infty e^{-s} s^{-(d-\beta)/4}\,ds\\
&\leq &C
\begin{cases}
1, & d-\beta>4,\\
1+\log(1/\alpha), & d-\beta=4,\\
\alpha^{(d-\beta)/4-1}, & d-\beta<4.
\end{cases}
\end{eqnarray*}
Combining these estimates with Lemma \ref{lem:theta}, it holds that
\begin{eqnarray*}
\int_0^\delta \int_E  R_\alpha^\mu f\, R_\alpha f\,d\mu_-\,d\widetilde{\nu}(\alpha) &\leq & C \int_0^\delta \alpha^{-1}\int_E  |R_\alpha f|^2\,V\,dx\,\alpha^{\beta/2}\,d\alpha \\
&\leq & C
\begin{cases}
\delta^{\beta/2}, & d\geq 6 \text{\ or\ }d=5, 0<\beta<1,\\
\delta^{1/2}(1+\log(1/\delta))^2, & d=5, \beta=1,\\
\delta^{1/2}, & d=5, 1<\beta<2.
\end{cases}
\end{eqnarray*}
Hence {\rm (UC)} holds. By Theorem \ref{thm:subcri_cri}, \(-(-\frac{1}{2}\Delta - \lambda(Vm)V)^{\beta/2}\) is critical for \(0<\beta \leq 2\) and \(d\geq 5\).
\end{example}

\appendix
\section{Dirichlet form}\label{append}
For the convenience, we summarize basic properties of closed forms and Dirichlet forms. See \cite{CF12,FOT11,O13} for more details.

Let \(E\) be a locally compact separable metric space and \(m\) be a positive Radon measure with \(\text{ supp}(m)=E\). The state space \(E\) is equipped with the Borel \(\sigma\)-algebra \(\mathcal{B}(E)\). The inner product in \(L^2(E;m)\) is denoted by \(\langle \cdot, \cdot \rangle_{m}\) and the \(L^2\)-norm is denoted by \(||\cdot ||_{L^2(m)}\). We take an isolated point \(\partial \not \in E\) called the cemetery point, and set \(E_\partial:=E\cup \{\partial\}\) equipped with \(\mathcal{B}(E_\partial):=\mathcal{B}(E) \cup \{B \cup \{\partial\} : B \in \mathcal{B}(E)\}\). 

\begin{definition}[closed form and Dirichlet form]
For a dense linear subspace \(\mathcal{D}(\mathcal{E})\) of \(L^2(E;m)\) and  a non-negative definite symmetric bilinear form \(\mathcal{E}\) on \(\mathcal{D}(\mathcal{E})\times \mathcal{D}(\mathcal{E})\), we call \(\mathcal{D}(\mathcal{E})\) the domain of \(\mathcal{E}\). A non-negative definite symmetric bilinear form \((\mathcal{E}, \mathcal{D}(\mathcal{E}))\) is called a \textit{closed form} on \(L^2(E;m)\) if \(\mathcal{E}\) is complete with respect to the norm induced by \(\mathcal{E}_1\), where \(\mathcal{E}_{\alpha}(f,g):=\mathcal{E}(f,g)+\alpha \langle f,g\rangle_m\) for \(\alpha>0\).

A closed form \((\mathcal{E}, \mathcal{D}(\mathcal{E}))\) is called a \textit{Dirichlet form} on \(L^2(E;m)\) if it is Markovian, that is, for any \(f\in \mathcal{D}(\mathcal{E})\), it holds that \(g:=(0\vee f)\wedge 1 \in \mathcal{D}(\mathcal{E})\) and \(\mathcal{E}(g,g)\leq \mathcal{E}(f,f)\).
\end{definition}
We note that \((\mathcal{D}(\mathcal{E}), \mathcal{E}_\alpha)\) is a Hilbert space for any closed form \((\mathcal{E}, \mathcal{D}(\mathcal{E}))\) and \(\alpha>0\).

For a strongly continuous contraction semigroup \(\{T_t\}_{t> 0}\) on \(L^2(E;m)\), the pair
\begin{eqnarray*}
\begin{cases}
\mathcal{D}(\mathcal{E}):=\{f \in L^2(E;m)\mid \lim_{t\searrow 0} \frac{1}{t}\langle f-T_tf, f \rangle_m <\infty\}\\
\mathcal{E}(f,g):=\lim_{t\searrow 0} \frac{1}{t}\langle f-T_tf, g \rangle_m\ \ \text{for\ }f,g\in \mathcal{D}(\mathcal{E})
\end{cases}
\end{eqnarray*}
is a closed form. Conversely, for a closed form \((\mathcal{E}, \mathcal{D}(\mathcal{E}))\) on \(L^2(E;m)\), there exists a strongly continuous contraction semigroup \(\{T_t\}_{t> 0}\) on \(L^2(E;m)\) such that \(G_{\alpha}f:= \int_0^{\infty} e^{-\alpha t} T_tf(x) dt \in \mathcal{D}(\mathcal{E})\) and \(\mathcal{E}_{\alpha}(G_{\alpha}f,g)=\langle f, g\rangle_m\) for \(\alpha >0\), \(f\in L^2(E;m)\) and \(g\in \mathcal{D}(\mathcal{E})\).

Furthermore, the generator of a strongly continuous contraction semigroup is a non-positive definite self-adjoint operator on \(L^2(E;m)\) and, for a non-negative definite self-adjoint operator \(-\mathcal{L}\) on \(L^2(E;m)\), \(\{T_t:=e^{\mathcal{L}t}\}_t\) is a strongly continuous contraction semigroup on \(L^2(E;m)\). In this case, the corresponding closed form is represented by \(\mathcal{D}(\mathcal{E})=\mathcal{D}(\sqrt{-\mathcal{L}})\) and \(\mathcal{E}(f,g)=\langle \sqrt{-\mathcal{L}}f, \sqrt{-\mathcal{L}}g \rangle_m = \langle -\mathcal{L}f, g \rangle_m\).

\begin{definition}[Markov process]
A quadruplet \(X=(\Omega, \mathcal{M}, \{X_t\}_{t\geq 0}, \{\mathbb{P}_x\}_{x\in E_\partial})\) is a Markov process on \(E\) if the following conditions hold.
\begin{enumerate}
\item[(M1)] For each \(x\in E_\partial\), \((\Omega, \mathcal{M}, \{X_t\}_{t\geq 0}, \mathbb{P}_x)\) is a stochastic process on \(E_\partial\), that is, \(\Omega, \mathcal{M}, \mathbb{P}_x\) is a probability space and \(X_t : \Omega \to E_\partial\) is a measurable map for each \(t\geq 0\).
\item[(M2)] For each \(t\geq 0\) and \(B\in \mathcal{B}(E_\partial)\), a map \(E_\partial \ni x \mapsto \mathbb{P}_x(X_t\in B) \in \mathbb{R}\) is measurable.
\item[(M3)] There exists a family of increasing sub \(\sigma\)-fields \(\{\mathcal{M}_t\}_{t\geq 0}\) of \(\mathcal{M}\) such that  \(X_t\) is \(\mathcal{M}_t\) measurable for each \(t\geq 0\), and \(\mathbb{P}_x(X_{s+t} \in B | \mathcal{M}_t) = \mathbb{P}_{X_t}(X_s\in B)\), \(\mathbb{P}_x\)-almost surely for any \(x\in E_\partial, s,t \geq 0\) and \(B\in \mathcal{B}(E_\partial)\). Here \(\mathbb{P}_x(\cdot|\cdot)\) is a conditional probability.
\item[(M4)] It holds that \(\mathbb{P}_x(X_0=x)=1\) for any \(x\in E_\partial\), and \(\mathbb{P}_\partial(X_t=\partial)=1\) for any \(t\geq 0\).
\end{enumerate}
\end{definition}
The condition (M3) is called the Markov property, and \(\{\mathcal{M}_t\}_{t\geq 0}\) in (M3) is called an admissible filtration. We call \(\sigma :\Omega \to [0,\infty]\) a stopping time if \(\{\sigma \leq t\} \in \mathcal{M}_t\) for each \(t\geq 0\). For a Markov process \(X\), we set \(\zeta(\omega):= \inf\{t \geq 0 : X_t(\omega)=\partial\}\), which is called a lifetime of \(X\).

\begin{definition}[Hunt process]
A Markov process \(X=(\Omega, \mathcal{M}, \{X_t\}_{t\geq 0}, \{\mathbb{P}_x\}_{x\in E_\partial})\) on \(E\) is called a Hunt process if the following conditions hold.
\begin{enumerate}
\item[(H1)] For \(t\geq \zeta(\omega)\), \(X_t(\omega) =\partial\) holds. For each \(t\geq 0\), there exists a map \(\theta_t :\Omega \to \Omega\) such that \(X_s\circ \theta_t = X_{s+t}\) for any \(s\geq 0\). For each \(\omega \in \Omega\), \(X_\cdot(\omega)\) is right continuous on \([0,\infty)\) and has left limits on \((0,\infty)\) in \(E_\partial\).
\item[(H2)]  An admissible filtration \(\{\mathcal{M}_t\}_{t\geq 0}\) in (M3) satisfies \(\bigcap_{s>t}\mathcal{M}_s = \mathcal{M}_t\) for each \(t\geq 0\) and, \(\mathbb{P}_\mu(X_{s+\sigma} \in B | \mathcal{M}_\sigma) = \mathbb{P}_{X_\sigma}(X_s\in B)\) holds \(\mathbb{P}_\mu\)-almost surely for each stopping time \(\sigma\), any \(s\geq 0\), \(B\in \mathcal{B}(E_\partial)\) and any probability measure \(\mu\) on \(E_\partial\).
\item[(H3)] For any increasing stopping times \(\{\sigma_n\}_n\) with \(\sigma := \lim_{n\to \infty} \sigma_n\), it holds that 
\(\mathbb{P}_\mu(\lim_{n\to \infty} X_{\sigma_m} = X_\sigma, \sigma<\infty) = \mathbb{P}_\mu(\sigma < \infty)\) for any probability measure \(\mu\) on \(E_\partial\).
\end{enumerate}
\end{definition}
The condition (H2) is called the strong Markov property, and a Markov process enjoying (H2) is called a strong Markov process. The condition (H3) is called quasi-left-continuity on \((0,\infty)\).

\begin{definition}[regular Dirichlet form]
A Dirichlet form \((\mathcal{E},\mathcal{D}(\mathcal{E}))\) is called \textit{regular} if \(\mathcal{D}(\mathcal{E}) \cap C_c(E)\) is dense in \((\mathcal{D}(\mathcal{E}), \mathcal{E}_1)\) and  dense in \((C_c(E), |\cdot||_{\infty})\), where \(\|\cdot\|_\infty\) is the essential supremum with respect to \(m\).
\end{definition}

\begin{definition}[nest, polar set, quasi-continuous, \(m\)-inessential set]\label{quasi_notation}
\ \vspace{-3mm}\\
\begin{enumerate}
\item An increasing sequence of closed sets \(\{F_k\}_{k\geq 1}\) of \(E\) is called a \textit{nest} if \(\cup_{k\geq 1} \{f\in \mathcal{D}(\mathcal{E}) : f=0 \ m\text{-a.e.\  on\  }E \setminus F_k\}\) is \(\mathcal{E}_1\)-dense in \(\mathcal{D}(\mathcal{E})\).
\item  \(N\subset E\) is \textit{\(\mathcal{E}\)-polar} if there exists a nest \(\{F_k\}_{k\geq 1}\) such that \(N\subset \cap_{k\geq 1} (E\setminus F_k)\).
\item A statement depending on \(x\in A \subset E\) holds \(\mathcal{E}\)-\textit{quasi-everywhere} (q.e. in abbreviation) on \(A\) if there exists an \(\mathcal{E}\)-polar set \(N\subset A\) such that the statement holds for \(x\in A\setminus N\).
\item A function \(f\) is \textit{quasi-continuous} if there exists a nest \(\{F_k\}_{k\geq 1}\) such that the restriction of \(f\) to \(F_k\) is finite and continuous on \(F_k\) for each \(k\geq 1\).
\item  A subset \(B\subset E\) is a \textit{nearly Borel set} if, for any probability measure \(\mu\) on \(E\cup \{\partial\}\), there exist Borel sets \(B_1, B_2\) such that \(B_1\subset B \subset B_2\) and \(\mathbb{P}_{\mu}(X_t \in B_2\setminus B_1 \text{\ for\ some\ }t\geq 0)=0\).
\item A subset \(N\subset E\) is \(m\)\textit{-inessential} if \(N\) is an \(m\)-negligible nearly Borel set such that \(\mathbb{P}_x(\sigma_{N}<\infty)=0\) for \(x\in E\setminus N,\) where \(\sigma_{N}:=\inf\{t>0;X_t\in N\}\) is the first hitting time to \(N\).
\end{enumerate}
\end{definition}

We remark that an \(\mathcal{E}\)-polar set is \(m\)-negligible.

For a regular Dirichlet form \((\mathcal{E},\mathcal{D}(\mathcal{E}))\) on \(L^2(E;m)\), by Fukushima's theorem, there exists an \(m\)-symmetric Hunt process \(X=(\Omega, \mathcal{M}, \{X_t\}_{t\geq 0}, \{\mathbb{P}_x\}_{x\in E_\partial})\) on \(E\) associated with \((\mathcal{E},\mathcal{D}(\mathcal{E}))\). Here, \(X\) is \(m\)-symmetric if it holds that
\[\int_E P_tf(x) g(x)\,dm(x) = \int_E f(x) P_tg(x)\,dm(x)\]
for any non-negative Borel measurable functions \(f,g\), where 
\[P_tf(x):=\mathbb{E}_x[f(X_t)]:=\int_\Omega f(X_t(\omega))\, d\mathbb{P}_x(\omega).\]
For \(f\in L^2(m)\) and \(\alpha>0\), we define resolvents \(\{R_{\alpha}\}_{\alpha>0}\) by
\[R_{\alpha}f(x):=\int_0^{\infty}e^{-\alpha t}P_tf(x)dt\]
 We remark that \(P_tf\) (resp. \(R_{\alpha}f\)) is a quasi-continuous version of \(T_tf\) (resp. \(G_{\alpha}f\)) for \(f\in \mathcal{D}(\mathcal{E})=\mathcal{D}(\sqrt{-\mathcal{L}})\), that is, for each \(t\geq 0\) and \(\alpha>0\), \(P_tf = T_tf \) and \(R_\alpha f = G_\alpha f\) \(m\)-almost everywhere.

\begin{example}
Let \(E:=\mathbb{R}^d\), \(m\) be a Lebesgue measure, and \(\mathcal{L}:=\Delta\). Then, the corresponding closed form \((\mathcal{E}, \mathcal{D}(\mathcal{E}))\) on \(L^2(\mathbb{R}^d; dx)\) is \(\mathcal{D}(\mathcal{E})=H^1(\mathbb{R}^d)\), the \(1\)-order Sobolev space, and \(\mathcal{E}(f,g)=\int \nabla f \, \nabla g\, dx\). In this case, \((\mathcal{E}, \mathcal{D}(\mathcal{E}))\) is a regular Dirichlet form on \(L^2(\mathbb{R}^d; dx)\) and its corresponding process is Brownian motion on \(\mathbb{R}^d\).
\end{example}

We define an \textit{extended Dirichlet space} \(\mathcal{D}_e(\mathcal{E})\) by the set of \(m\)-equivalence classes of all \(m\)-measurable functions \(f\) on \(E\) such that \(|f|<\infty, m\)-almost everywhere and having an \(\mathcal{E}\)-Cauchy sequence
\(\{f_n\}_{n\geq 1}\subset \mathcal{D}(\mathcal{E})\) converging to \(f\) \(m\)-almost everywhere on \(E\). We can define \(\mathcal{E}(f,f)\) for \(f\in \mathcal{D}_e(\mathcal{E})\) by \(\mathcal{E}(f,f):= \lim_{n\to \infty}\mathcal{E}(f_n,f_n)\) for the above sequence \(\{f_n\}_n \subset \mathcal{D}(\mathcal{E})\). This definition is independent of the choice of an approximation sequence \(\{f_n\}_n\) of \(f\in \mathcal{D}_e(\mathcal{E})\). We remark that any function belonging to \(\mathcal{D}(\mathcal{E})\) has a quasi-continuous version, so without loss of generality, we may treat all functions in \(\mathcal{D}_e(\mathcal{E})\) as quasi-continuous functions.

\begin{definition}[Transience and Recurrence]
A regular Dirichlet form \((\mathcal{E}, \mathcal{D}(\mathcal{E}))\) on \(L^2(E;m)\) is {\it transient} if \((\mathcal{D}_e(\mathcal{E}), \mathcal{E})\) is a real Hilbert space. A regular Dirichlet form \((\mathcal{E}, \mathcal{D}(\mathcal{E}))\) on \(L^2(E;m)\) is {\it recurrent} if \(1\in \mathcal{D}_e(\mathcal{E})\) and \(\mathcal{E}(1,1)=0.\)
\end{definition}

\begin{definition}[Irreducibility]
For a strongly continuous semigroup \(\{T_t\}_t\), a set \(A\) is \(\{T_t\}\)-invariant if \(T_t(1_A f)=1_A T_tf\) for any \(f\in L^2(E;m)\) and \(t>0\). A semigroup \(\{T_t\}_t\) is irreducible if, for any \(\{T_t\}\)-invariant set \(A\), either \(m(A)=0\) or \(m(A^c)=0\) holds. 
\end{definition}

We also say a closed form \((\mathcal{E}, \mathcal{D}(\mathcal{E}))\) is irreducible if its corresponding semigroup is irreducible.

\begin{proposition}[{cf. \cite[Proposition 2.1.3]{CF12}}]
An irreducible Dirichlet form is either transient or recurrent.
\end{proposition}

\begin{definition}[Smooth measure]
A positive Borel measure \(\mu\) on \(E\) is a \textit{smooth measure} if \(\mu\) charges no \(\mathcal{E}\)-polar set and, there exists a nest \(\{F_k\}_k\) such that \(\mu(F_k)<\infty\) for every \(k\geq 1.\) Denote by \(\mathcal{S}\) the family of all smooth measures.
\end{definition}

Let \(\mathcal{M}_{\infty}\) be the smallest \(\sigma\)-field including \(\{\mathcal{M}_t\}_t\), which is an admissible filtration of \(X\).
\begin{definition}[PCAF, positive continuous additive functional]\label{DefPCAF}
A \([-\infty,\infty]\)-valued stochastic process \(A=\{A_t\}_{t\geq 0}\) is called an \textit{additive functional} of \(X\) if there exist \(\Lambda \in \mathcal{M}_{\infty}\) and an \(m\)-inessential set \(N\subset E\) such that \(\mathbb{P}_x(\Lambda)=1\) for \(x\in E\setminus N\) and \(\theta_t \Lambda \subset \Lambda\) for any \(t>0\), and the following conditions hold.
\begin{enumerate}
\item[(A.1)] For each \(t\geq 0,\) \(A_t|_{\Lambda}\) is \(\mathcal{M}_t|_{\Lambda}\)-measurable.
\item[(A.2)] For any \(\omega \in \Lambda\), \(A_{\cdot}(\omega)\) is right continuous on \([0,\infty)\) and has left limits on \((0, \zeta(\omega))\), \(A_0(\omega)=0,\) \(|A_{t}(\omega)|<\infty\) for \(t<\zeta(\omega)\) and \(A_t(\omega)=A_{\zeta(\omega)}(\omega)\) for \(t\geq \zeta (\omega)\). Moreover the equation
\[A_{t+s}(\omega)=A_t(\omega)+A_s(\theta_t\omega)\ \ \text{for\ every\ }t,s\geq 0,\]
is satisfied.
\end{enumerate}

An additive functional \(A\) is called a \textit{positive continuous additive functional} (PCAF in abbreviation) if \(A\) is a \([0,\infty]\)-valued continuous process, and denote by \({\bf A}_c^+\) the family of all PCAFs.
\end{definition}
The set \(\Lambda\) appearing in Definition \ref{DefPCAF} is called the defining set of \(A\). A PCAF \(A\) is called a PCAF in the strict sense if \(N\) appearing in Definition \ref{DefPCAF} is empty. PCAFs \(A\) and \(B\) are called \textit{\(m\)-equivalent} if \(\int_E\mathbb{P}_x(A_t\neq B_t)dm(x)=0\) for any \(t>0\). An \(m\)-equivalence is equivalent to the existence of a common defining set \(\Lambda\) and a common \(m\)-inessential set \(N\) such that \(A_t(\omega)=B_t(\omega)\) for any \(t \geq 0\) and \(\omega \in \Lambda.\) 

PCAFs and smooth measures are in one-to-one correspondence in the following sense. This correspondence is called \textit{the Revuz correspondence}. Therefore, a smooth measure is also called the \textit{Revuz measure}. See \cite[Theorem 4.1.1]{CF12} for details. For a measure \(\rho,\) we set \(\mathbb{E}_\rho[\cdot]:=\int_E \mathbb{E}_x[\cdot] \,d\rho(x).\)
\begin{theorem}[The Revuz correspondence]
\((i)\) For a PCAF \(A\), there exists a unique smooth measure \(\mu\) such that
\begin{equation}
\int_E fd\mu = \lim_{t \to 0}\frac{1}{t} \mathbb{E}_m\left[ \int_0^t f(X_s) dA_s  \right] \label{eq:AppPCAF-1}
\end{equation}
for any non-negative Borel function \(f\) on \(E.\)\\
\((ii)\) For any smooth measure \(\mu\), there exists a PCAF \(A\) satisfying \((\ref{eq:AppPCAF-1})\) up to the \(m\)-equivalence.
\end{theorem}

For example, for a bounded positive Borel function \(f\), we set
\[A_t:=\int_0^t f(X_s)ds,\]
then \(A:=\{A_t\}_t\) is a PCAF and \(A\) corresponds to a smooth measure \(fdm\). As another example, when a capacity of \(x_{0}\in E\) is positive, the local time \(L^{x_0}\) is a PCAF and its corresponding smooth measure is a Dirac measure \(\delta_{x_0}\). Denote by \(A^\mu\) a PCAF corresponding to \(\mu\).

\bigskip
\noindent
{\large \bf Acknowledgments.}

\noindent
The author thanks to Dr. Takahiro Mori for fruitful discussions. This work was supported by JSPS KAKENHI Grant Number 25K17270.

\end{document}